\documentclass[graybox]{SNmult}

\usepackage{stmaryrd}          %
\SetSymbolFont{stmry}{bold}{U}{stmry}{m}{n} %
\makeatletter\let\@font@warning\@gobble\makeatother
\usepackage{amsmath,amssymb}   %
\usepackage{mathtools}         %
\mathtoolsset{centercolon}
\usepackage{microtype}         %
\usepackage{esint}             %
\usepackage{xfrac}             %
\usepackage{graphicx}          %
\usepackage{verbatim}

\DeclarePairedDelimiter{\abs}{\lvert}{\rvert}

\DeclarePairedDelimiter{\bra}{(}{)}

\DeclarePairedDelimiter{\set}{\{}{\}}

\newcommand{\Z}{\mathbb{Z}}
\newcommand{\R}{\mathbb{R}}
\newcommand{\Rp}{\mathbb{R}_{>0}}
\newcommand{\Rnn}{\mathbb{R}_{\ge0}}
\newcommand{\dd}{\,\mathrm{d}}
\newcommand{\dn}{\mathrm{d}}

\newcommand{\ent}{\mathcal{H}}
\DeclareMathOperator{\Ent}{Ent}
\newcommand{\logmean}{\boldsymbol{\Lambda}}
\newcommand{\geom}{\boldsymbol{\Gamma}}
\newcommand{\hell}{\mathbb{H}}
\newcommand{\wass}{\mathbb{W}}
\newcommand{\metric}{{\mathbb W\!\mathbb D}}

\newcommand{\mob}{\mathbf{m}}
\newcommand{\Mob}{\mathbf{M}}

\newcommand{\dff}{\operatorname{D\!}{}}

\newcommand{\dprb}{{\mathcal P}_2(\delta\mathbb{Z}^*)}
\newcommand{\dprbK}{{\mathcal P}_2(\delta\llbracket-K,K\rrbracket^*)}
\newcommand{\grid}{\delta\mathbb{Z}}
\newcommand{\kgrid}{\llbracket -K,K\rrbracket}

\newcommand{\hk}{{\mathbb H\!\mathbb K}}
\newcommand{\geoLap}{\Delta^{\mathrm{geo}}}

\newcommand{\GS}{Giuseppe Savar\'{e}}

\newtheorem{manualresultinner}{Result}
\newenvironment{manualresult}[1]{%
  \renewcommand\themanualresultinner{#1}%
  \manualresultinner
}{\endmanualresultinner}

\title*{The confined discrete DLSS equation:
        gradient flow structure and long-time asymptotics}
\titlerunning{The confined discrete DLSS equation}

\author{Daniel Matthes\orcidID{0000-0001-6764-8301}
   \and Andr\'e Schlichting\orcidID{0000-0003-4140-491X}}
\authorrunning{D.\ Matthes and A.\ Schlichting}

\institute{%
  Daniel Matthes
  \at Zentrum Mathematik, TU München, Boltzmannstrasse 3, D-85748 Garching, Germany.
  \email{matthes@ma.tum.de}
  \and
  Andr\'e Schlichting
  \at Institute for Applied Analysis, Ulm University, Helmholtzstrasse 18, D-89081 Ulm, Germany.
  \email{andre.schlichting@uni-ulm.de}
}

\usepackage[hidelinks,bookmarks=false]{hyperref}
\hypersetup{
  pdfauthor={Daniel Matthes, André Schlichting},
  pdftitle={The confined discrete DLSS equation:
        gradient flow structure and long-time asymptotics},
}
\makeatletter
\providecommand*{\toclevel@title}{0}
\providecommand*{\toclevel@author}{0}
\makeatother

\begin{document}

\maketitle

\abstract{%
We propose and analyse a structure-preserving discretization of the
fourth-order DLSS equation in dynamically rescaled variables leading to a linear confinement
on the one-dimensional real line.
The scheme is  obtained as a discrete analogue of the
Carrillo--Toscani splitting. It allows for a
new interpretation of the DLSS equation as gradient flow of the
relative entropy with respect to a combined diffusive and transport metric,
extending the structure of the discretization introduced in our earlier
work with Giuseppe Savar\'{e} for the unconfined equation. \newline
We prove global existence of positive solutions and, as the main
analytical result, exponential decay of the relative entropy towards
the discrete Gaussian equilibrium. The decay rate is obtained via a
discrete modified logarithmic Sobolev inequality.
The resulting rate is consistent with the continuum rate and converges to it from below in the vanishing mesh-size limit. 
Numerical results illustrate and confirm our findings.
}

\newpage 

\section{Introduction}
The fourth-order DLSS equation
\begin{align}
  \label{eq:DLSS0}
  \partial_t\rho = -\partial_{xx}\bigl(\rho\,\partial_{xx}\log\rho\bigr)
\end{align}
was first derived in the context of interface fluctuations in certain spin systems~\cite{DLSS1,DLSS2},
but also appears in the contexts of diffusive quantum hydrodynamics~\cite{DMR}
and, more recently, that of chemical reactions~\cite{MSS}.
Its well-posedness on bounded convex domains in any dimension was revisited by \GS\ and the authors~\cite{MatthesSavareSchlichting2026}.
In the square-root variable $u=\sqrt\rho$, the DLSS operator has a unique maximal monotone extension in $L^2$,
which generates a Hellinger-contractive semigroup.
In two papers~\cite{GST,MMS}, \GS\ has studied --- in a much more general context --- the well-posedness and long-time asymptotics of non-negative unit-mass solutions
of the following variant of~\eqref{eq:DLSS0} with linear confinement on the real line $\R$:
\begin{align}
  \label{eq:DLSS1}
  \partial_t\rho = -\partial_{xx}\bigl(\rho\,\partial_{xx}\log\rho\bigr) + \lambda^2\partial_x(x\,\rho).
\end{align}
The additional confinement can be interpreted as passage to the self-similar variables in~\eqref{eq:DLSS0}.
It was shown in~\cite{MMS} that solutions converge at the uniform exponential rate $2\lambda^2$
to the Gaussian equilibrium profile
\begin{align}
  \label{eq:Gauss}
    \bar\rho(x) = \frac{1}{Z} \exp\biggl(-\frac{\lambda}{2} x^2\biggr) \quad\text{with}\quad Z:= \sqrt{2\pi/\lambda}
\end{align}
in various metrics, implying self-similar asymptotics for the unconfined equation~\eqref{eq:DLSS0}.
In this paper, we propose a spatial discretization of~\eqref{eq:DLSS1} that replicates that equilibration behaviour.

To motivate the rather particular choice for our discretization,
we briefly recall a fundamental observation from~\cite{CarTos}
that is the key to the proof for equilibration in~\cite{MMS}.
Namely, one can rewrite~\eqref{eq:DLSS1}
in the form of a (combined first- and second-order) continuity equation as follows,
\begin{align}
  \label{eq:DLSSflux}
  \partial_t\rho = \partial_{xx}I[\rho] - \lambda\partial_xJ[\rho]
\end{align}
with second- and first-order fluxes $I$ and $J$ given, respectively, by
\begin{align}
  \label{eq:DLSS-IJ}
  I[\rho] = -\rho\,\partial_{xx}\log\frac{\rho}{\bar\rho},
  \quad
  J[\rho] = -\rho\,\partial_{x}\log\frac{\rho}{\bar\rho}.
\end{align}
To see that~\eqref{eq:DLSSflux} is indeed an equivalent reformulation of~\eqref{eq:DLSS1},
expand the fourth- and the second-order terms in~\eqref{eq:DLSSflux}
and observe that the arising second-order diffusion terms cancel.
The formal argument by which~\eqref{eq:DLSSflux}\&\eqref{eq:DLSS-IJ}
provides a quantitative control on entropy decay
--- and eventually results in the asymptotic estimate~\eqref{eq:Hdecay} below ---
is presented in Section~\ref{sct:cartos}.

Our spatial discretization builds directly on~\eqref{eq:DLSSflux}.
We discretize the continuous density $\rho$
by one that is piecewise constant on the cells of the infinite grid $\delta\Z$ with uniform spacing $\delta>0$.
With the notation introduced in Section~\ref{sct:discrete} below,
the discretized evolution takes the form
\begin{subequations}\label{eq:scheme}
\begin{equation}\label{eq:dCE2}
  \dot{\rho} = \Delta I^\delta[\rho] - \lambda \dff^* J^\delta[\rho] \,.
\end{equation}
with discrete fluxes $I^\delta[\rho]$ and $J^\delta[\rho]$ given at half-integer and integer sites $\kappa\in \grid^*$ and $k\in \grid$,
respectively, by
\begin{align} %
  \label{eq:def:diffFlux}
  I_\kappa^\delta[\rho] &
                          =  - \frac{2}{\delta^2}\bigl( e^{\frac{\lambda}{2} \delta^2} \sqrt{\rho_{\kappa-1}\cdot \rho_{\kappa+1}} - \rho_{\kappa}\bigr) \\
  \label{eq:def:transpFlux}
  J_k^\delta[\rho] &
              = - \frac{1}{\delta}\bra*{ e^{\frac{\lambda}{2} \delta x_k} \rho_{k+1/2}- e^{-\frac{\lambda}{2} \delta x_k} \rho_{k-1/2}}.
\end{align}
\end{subequations}
The particular discretization of the diffusive flux $I$ is based on a weighted version of the \emph{geometric Laplacian} motivated by recent joint work~\cite{MRSS} of \GS\ with the authors of the contribution at hand. There, a highly structure-preserving discretization of the unconfined DLSS equation~\eqref{eq:DLSS0}
on the one-dimensional torus has been introduced.
The discretization of the first-order flux $J$ in~\eqref{eq:DLSS-IJ} is inspired by~\cite{Maa, Mielke2013}.
In Section~\ref{sct:consistency}, we provide a short formal calculation
that indicates consistency of the discretized fluxes $I^\delta$ and $J^\delta$ from~\eqref{eq:def:diffFlux}\&\eqref{eq:def:transpFlux}
with their continuous counterparts $I$ and $J$ in~\eqref{eq:DLSS-IJ}.

A related asymptotics-preserving spatial discretization
  of the DLSS equation with confinement potential (on the multi-dimensional cube)
  has been proposed by the first author with Maas in~\cite{MaaM}.
  That scheme is not constructed as a discrete realization of the splitting from \cite{CarTos},
  but rather on the grounds of an observation by Denzler and McCann from~\cite{DenMcC}
  about the convexity-like behaviour of the (totally non-displacement-convex) Fisher information
  along certain geodesics in Wasserstein space.
  Just as for the scheme~\eqref{eq:dCE2} under consideration here,
  the one from~\cite{MaaM} dissipates the discrete relative entropy at the appropriate exponential rate,
  and in addition even the discrete relative Fisher information.
  The drawback of the discretization from~\cite{MaaM}, however,
  is its very abstract definition, which apparently makes 
  rigorous proofs of global well-posedness or convergence difficult,
whereas such an analysis is feasible for~\eqref{eq:dCE2}
thanks to the results in~\cite{MRSS} obtained by the authors in collaboration with \GS.

We summarize the three main results of this paper. 
We emphasize that the results consider the setting of fixed $\lambda>0$ and $\delta>0$; for this reason the explicit dependence on $\lambda$ is omitted in the notation, whereas superscript $\delta$ is used to indicate quantities related to the discretization.
The first result concerns the variational structure of the proposed scheme.
\begin{svgraybox}
  \begin{manualresult}{A}[Gradient flow structure]\label{res:GF}
    There are mobility functions~$\Mob^\delta$ and~$\mob^\delta$, homogeneous of degree one in $\rho$,
    such that~\eqref{eq:scheme} can be written in gradient flow form
    \begin{align}
      \label{eq:dDLSSgf}
      \dot \rho = -\mathcal{K}^\delta(\rho)\bigl[\bigl(\ent_\lambda^\delta\bigr)'(\rho)\bigr],
    \end{align}
    with the positive semi-definite operator
    \begin{align}
      \label{eq:dOnsager} 
      \mathcal{K}^\delta(\rho)[\xi] = \Delta\bigl(\Mob^\delta[\rho] \,\Delta\xi\bigr) - \lambda \dff^*\bigl(\mob^\delta[\rho] \dff\xi\bigr) \quad\text{ for } \xi : \grid^* \to \R \,.
    \end{align}\vspace{-1.5\baselineskip}
  \end{manualresult}
\end{svgraybox}
The second result establishes global well-posedness on the infinite lattice;
below, $\dprb$ is the space of probability measures with finite second moment
on the $\delta$-spaced lattice.
\begin{svgraybox}
  \begin{manualresult}{B}[Global existence and positivity]\label{res:exists}
    For every non-negative initial datum $\rho^0\in\dprb$
    the scheme~\eqref{eq:scheme} admits a global solution $\rho^*:[0,\infty)\to\dprb$.
    Moreover, $\rho^*$ is strictly positive for all $t>0$.
  \end{manualresult}
\end{svgraybox}
The third and main analytical result is exponential convergence to equilibrium
with an explicit, $\delta$-uniform rate.
In the context of~\cite{MMS}, the main measure for proximity
of a given probability density $\rho$ to the equilibrium state $\bar\rho$, a Gaussian with variance $1/\lambda$, has been the relative entropy,
\begin{align}
  \label{eq:H}
  \ent_\lambda(\rho) = \int \bra*{ \rho \log\frac\rho{\bar\rho} - \rho + \bar\rho }\dd x,
\end{align}
which is non-negative and zero iff $\rho=\bar\rho$.
It was shown that solutions to the confined DLSS equation~\eqref{eq:DLSS1} satisfy
\begin{align}
  \label{eq:Hdecay}
  \ent_\lambda\bigl(\rho(t)\bigr) \le e^{-2\lambda^2t}\ent_\lambda\bigl(\rho(0)\bigr).
\end{align}
For our spatially discretized evolution~\eqref{eq:scheme},
we identify a steady state $\bar\rho^\delta\in\dprb$ in~\eqref{eq:dGauss}
and introduce a discrete version $\ent_\lambda^\delta$ of the relative entropy in~\eqref{eq:dH}.
The following is our analogue of~\eqref{eq:Hdecay}.
\begin{svgraybox}
  \begin{manualresult}{C}[Exponential convergence]\label{res:longtime}
    Let $\rho^*:[0,\infty)\to\dprb$ be the solution from Result~\ref{res:exists}.
    Then the relative entropy decays exponentially:
    \begin{equation}
      \label{eq:dHdecay}
      \ent_\lambda^\delta\bigl(\rho^*(t)\bigr)
      \;\leq\;
      e^{-2\lambda \alpha^\delta_\lambda t}\,
      \ent_\lambda^\delta(\rho^0)
      \qquad\text{for all }t\geq 0\,,
    \end{equation}
    where the explicit rate constant is given as
    \begin{equation}\label{eq:alphaLSI}
      \alpha^\delta_\lambda
      := \frac{2}{\delta^2}\bigl(1-e^{-\frac\lambda2\delta^2}\bigr)  = \lambda \biggl(1 - \frac{\lambda}{4} \delta^2 + O(\delta^4)\biggr)  \,.
    \end{equation}
    In particular, $\alpha^\delta_{\lambda}\to\lambda$ as $\delta\to 0$.
  \end{manualresult}
\end{svgraybox}
We see that the discrete rate equals the continuum rate $2\lambda^2$
established in~\cite{MMS} for the PDE~\eqref{eq:DLSS1}, up to errors of order $\delta^2$.
The proof proceeds via a discrete modified logarithmic Sobolev inequality,
obtained by identifying the transport dissipation with the Dirichlet form
of an explicitly computable birth--death Markov chain with stationary
measure $\bar\rho^\delta$ and applying the discrete Bakry--\'{E}mery criterion
of~\cite{CaputoPraPosta2009}.

\section{Analytic and numerical examples}\label{sec:AnaNumExamples}
\begin{figure}
	\centering
	\includegraphics[width=0.495\linewidth]{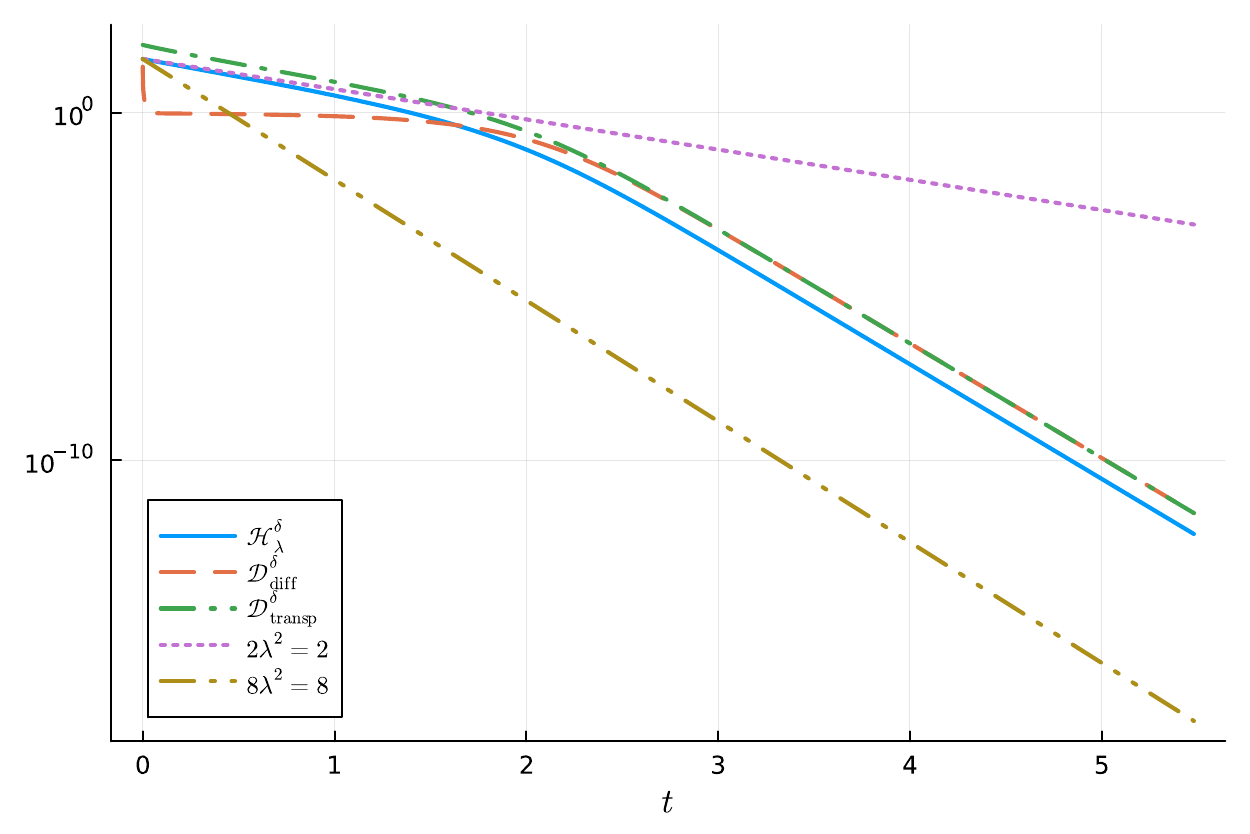}
	\hfill
	\includegraphics[width=0.495\linewidth]{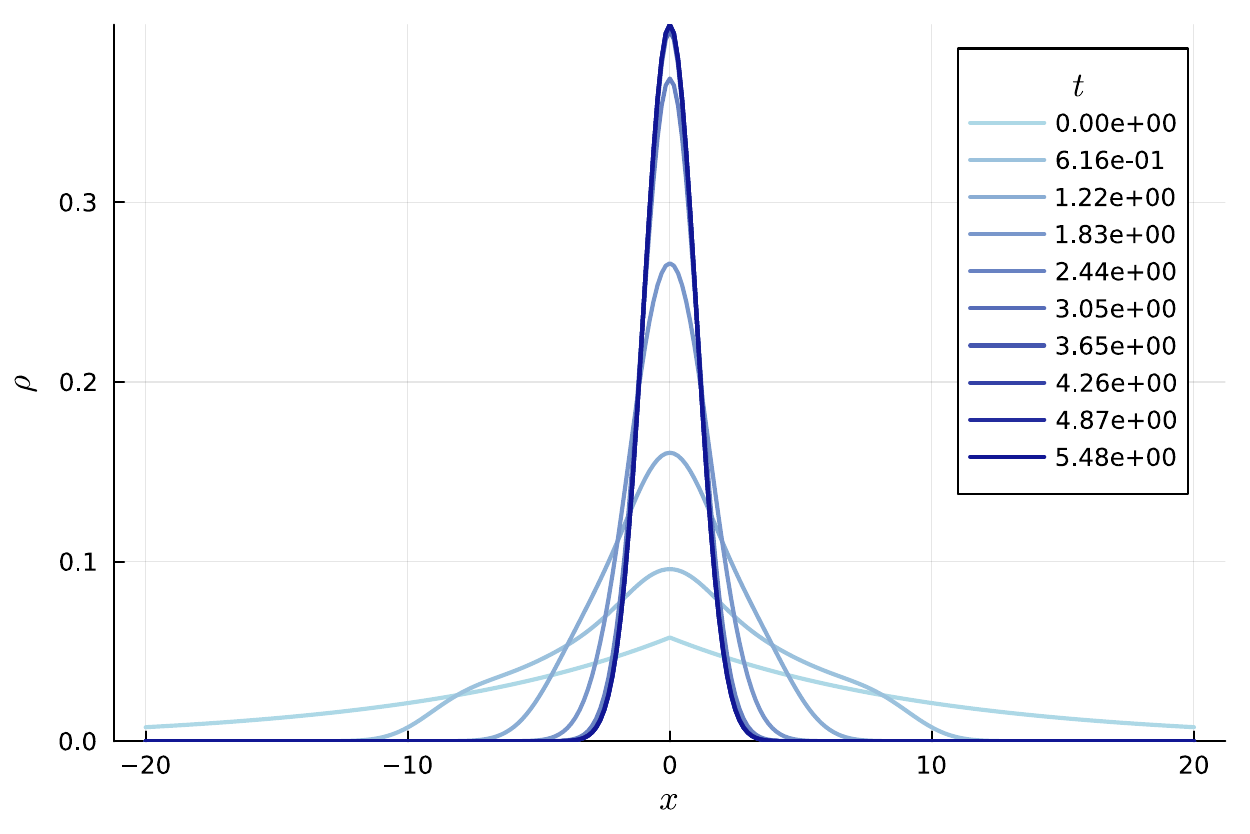}
	\caption{Numerical solution started from centred exponential profile $\rho_0 \propto e^{-\abs{x}/10}$ on domain $[-20,20]$. The convergence happens like in the centred Gaussian case with improved rate $8\lambda^2$.}\label{fig:ExpCentered}
\end{figure}

\begin{figure}
	\centering
	\includegraphics[width=0.495\linewidth]{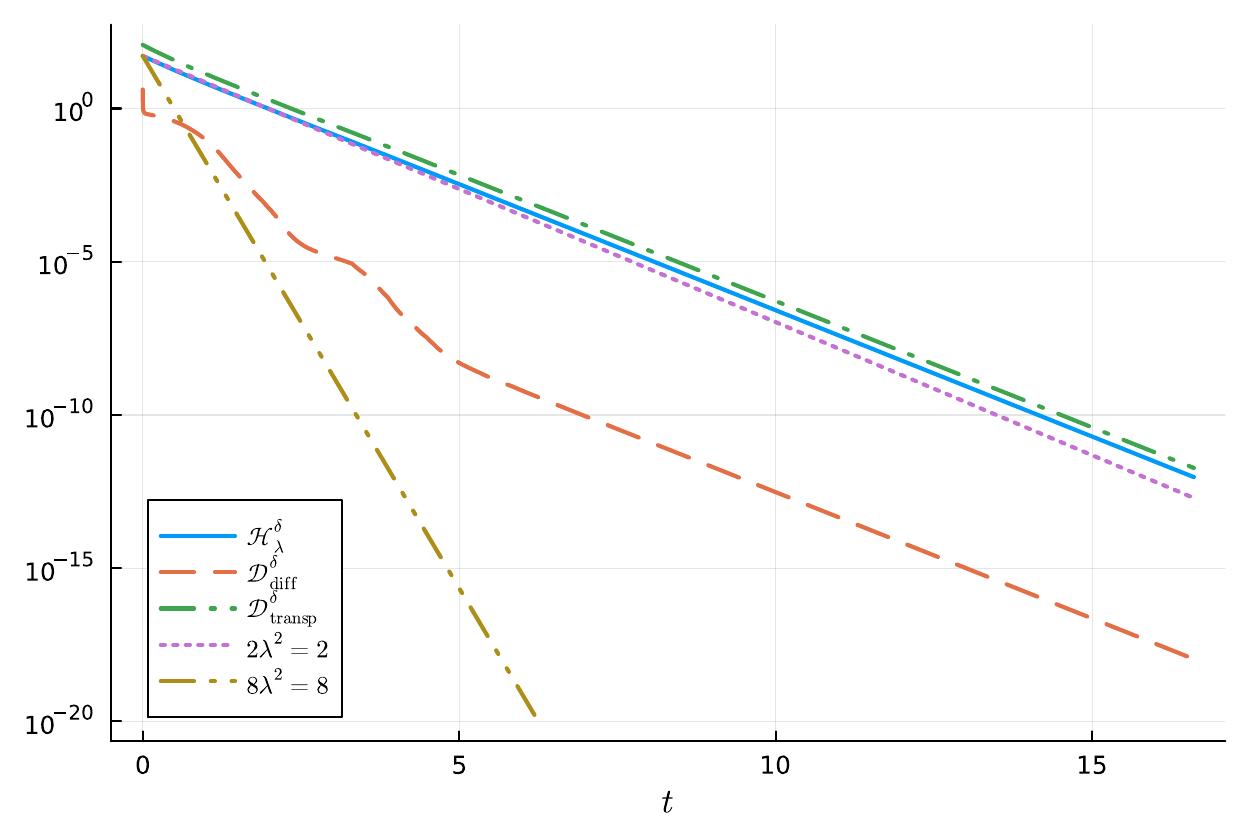}
	\hfill
	\includegraphics[width=0.495\linewidth]{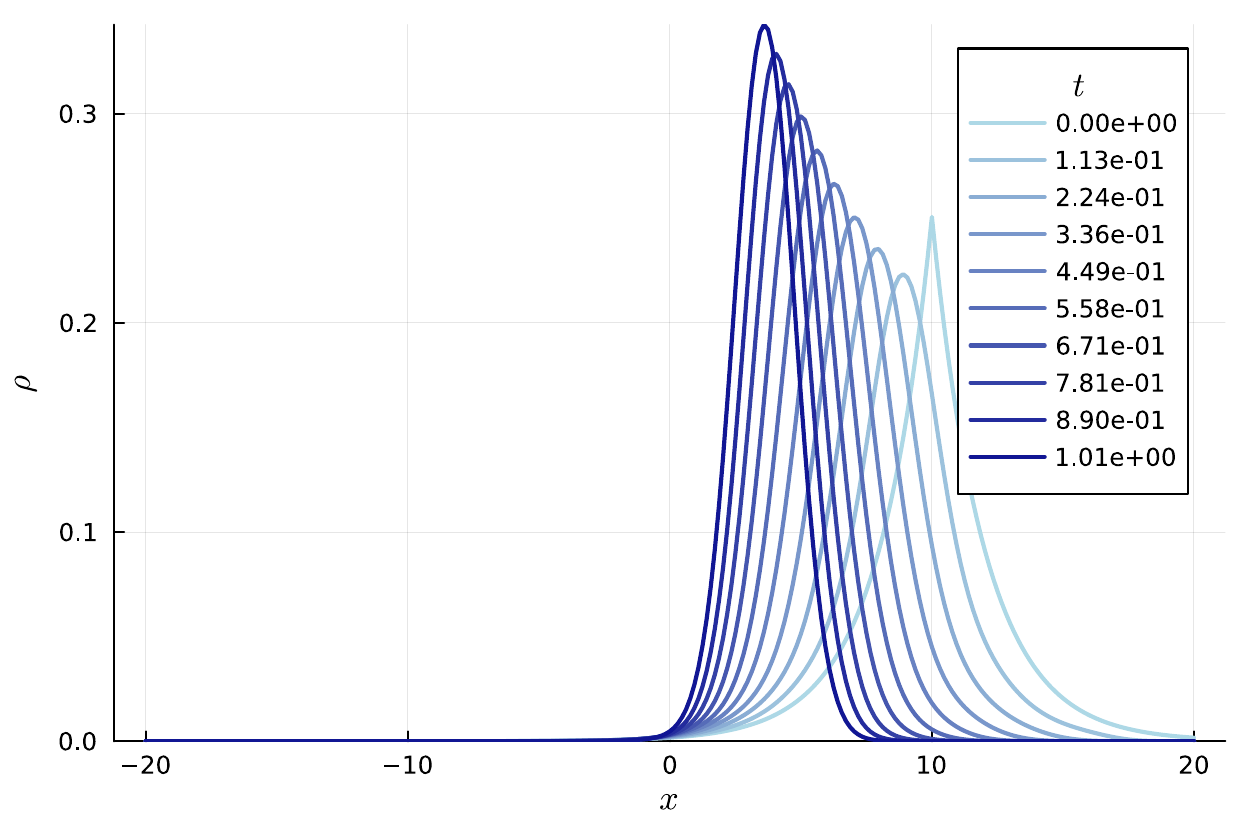}
	\caption{Numerical solution started from shifted exponential profile $\rho_0 \propto e^{-\abs{x-10}/2}$ on domain $[-20,20]$. Off-centred data attain the theoretical rate $2\lambda^2$.}\label{fig:ExpShifted}
\end{figure}

\begin{figure}
		\centering
	\includegraphics[width=0.495\linewidth]{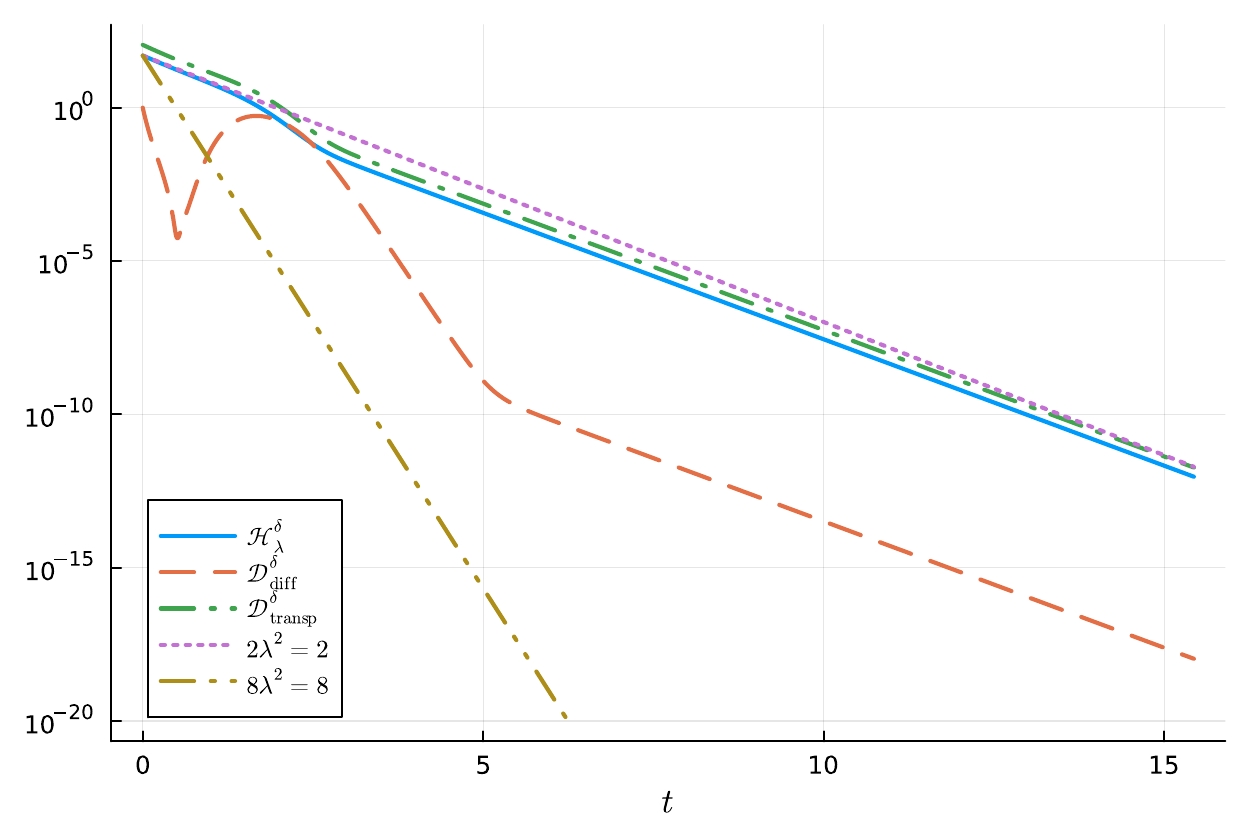}
	\hfill
	\includegraphics[width=0.495\linewidth]{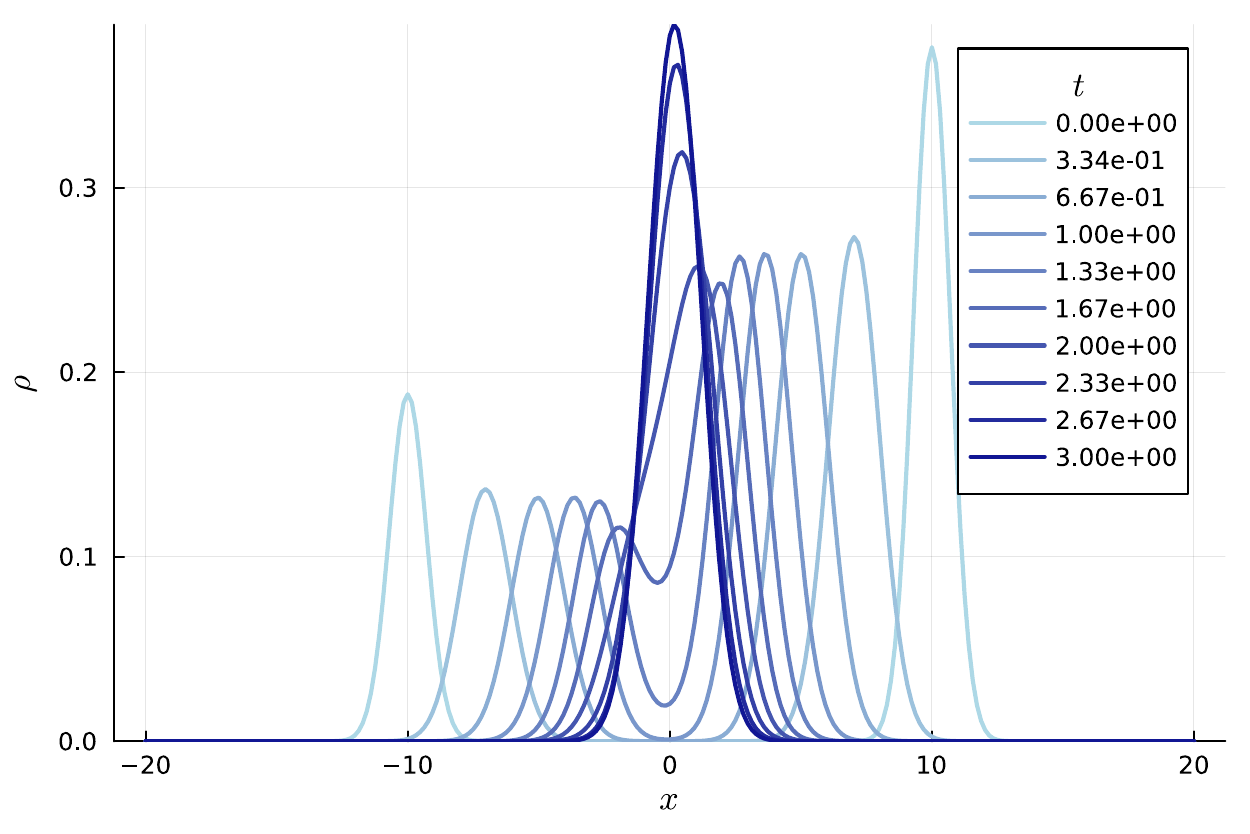}
	\caption{Numerical solution started from asymmetric sum-of-Gaussians profile $\rho_0 \propto 1.5 e^{-\abs{x-10}^2}+ e^{-\abs{x+10}^2}$ on domain $[-20,20]$. The diffusive part of the dissipation $\mathcal{D}_{\mathrm{diff}}$ (see~\eqref{eq:dissipation:diff} below) is not monotone along the evolution.}\label{fig:2GaussAsym}
\end{figure}

Before entering into details of the discretization and its analysis,
we would like to provide some intuition for the long-time asymptotics
of the confined DLSS equation~\eqref{eq:DLSS1} in this section.
Specifically, we consider the temporal behaviour of certain quantities along selected solutions $\rho(t)$.
The quantities of interest are the relative entropy $\ent_\lambda(\rho(t))$ from~\eqref{eq:H}
and its total dissipation along solutions to~\eqref{eq:DLSSflux},
\begin{align*}
  \mathcal D(\rho):=
  -\frac{\dd}{\dn t}\ent_\lambda(\rho)
  = -\int\partial_t\rho\log\frac{\rho}{\bar\rho}
  =  {\mathcal D}_\text{diff}(\rho) +  \lambda {\mathcal D}_\text{transp}(\rho),
\end{align*}
with the individual contributions by diffusion and by transport, i.e.,
\begin{align*}
  {\mathcal D}_\text{diff}(\rho) :=  \int \rho\left[\partial_{xx}\log\biggl(\frac{\rho}{\bar\rho}\biggr)\right]^2 \dd x
  \quad \text{and} \quad
  {\mathcal D}_\text{transp}(\rho) :=  \int \rho\left[\partial_{x}\log\biggl(\frac{\rho}{\bar\rho}\biggr)\right]^2 \dd x.
\end{align*}
As a benchmark for our numerical experiments,
we use a class of explicit solutions to the confined DLSS equation \eqref{eq:DLSS1},
given by scaled Gaussians as follows.
\begin{proposition}[Gaussian solutions]\label{prop:Gaussian}
  For mean $m(t)$ and variance $v(t)$ given by
  \begin{equation}\label{eq:GaussSol}
    m(t) = m_0\,e^{-\lambda^2 t},
    \qquad
    v(t) = \left[\frac{1}{\lambda^2}
      + \left(v_0^2 - \frac{1}{\lambda^2}\right)e^{-4\lambda^2 t}\right]^{1/2},
  \end{equation}
  the corresponding Gaussians
  \begin{equation}\label{eq:GaussSol:form}
    \rho(t,x) = \frac{1}{\sqrt{2\pi v(t)}}
    \exp\!\left(-\frac{(x-m(t))^2}{2v(t)}\right),
  \end{equation}
  satisfy~\eqref{eq:DLSS1}.
  Moreover, $m(t)\to 0$ and $v(t)\to 1/\lambda$ as $t\to\infty$,
  so $\rho(t)\to\bar\rho$.
\end{proposition}
\begin{proof}
  With the ansatz~\eqref{eq:GaussSol:form},
  one finds $\partial_{xx}\log\rho \equiv -1/v$, hence $I[\rho]= \rho(1/v-\lambda)$.
  The transport term and time derivative $\partial_t\rho$
  generate factors of the type $\rho(x-m)^\ell$ for $\ell=0,1,2$,
  and matching coefficients in the equation
  yields the decoupled ODEs
  \begin{equation*}
    \dot m = -\lambda^2 m
    \quad\text{with } m(0)=m_0
    \qquad\text{and}\qquad
    \dot v = \frac{2}{v} - 2\lambda^2 v
    \quad\text{with } v(0)=v_0 \,.
  \end{equation*}
  The equation for $m$ is linear. For $v$,
  the substitution $w = v^2$ gives $\dot w = 4 - 4\lambda^2 w$,
  which integrates to~\eqref{eq:GaussSol}.
\end{proof}
For the Gaussians $\rho(t)$ from~\eqref{eq:GaussSol:form},
the relative entropy to the steady state $\bar\rho$ has the explicit expression
\begin{equation}\label{eq:GaussEnt}
  \ent_\lambda\bigl(\rho(t)\bigr)
  = \frac{1}{2}\bigl[\lambda v(t) - 1 - \log(\lambda v(t))\bigr]
  + \frac{\lambda}{2}\,m(t)^2.
\end{equation}
Substituting~\eqref{eq:GaussSol} yields a fully explicit formula for the entropy decay.
There are several interesting regimes:
\begin{itemize}
\item \emph{Off-centre data, $m_0\neq 0$, $v_0 = 1/\lambda$}:
  The variance is constant $v(t)\equiv 1/\lambda$, so the entropy
  reduces to $\frac{\lambda}{2}m_0^2 e^{-2\lambda^2 t}$,
  decaying exactly at the general rate~$2\lambda^2$.
  In this case the diffusive second-order dissipation vanishes identically.		
\item \emph{Centred data, $m_0=0$, $v_0\neq 1/\lambda$}:
  The entropy is $\frac{1}{2}[\lambda v(t)-1-\log(\lambda v(t))]$.
  Near equilibrium, $v(t)-1/\lambda \sim C e^{-4\lambda^2 t}$,
  so by~\eqref{eq:GaussEnt} the entropy decays as $Ce^{-8\lambda^2 t}$,
  \emph{four times faster} than the general lower bound $2\lambda^2$.
  This improvement is entirely due to the diffusive dissipation term, which contributes
  a Fisher information of order $(1/v - \lambda)^2$ when $v\neq 1/\lambda$.  
\item \emph{General data, $m_0\neq 0$, $v_0\neq 1/\lambda$}:
  The entropy~\eqref{eq:GaussEnt} is a sum of the two terms above,
  decaying at rates $2\lambda^2$ and $8\lambda^2$ respectively.
  For large times the slower mean term dominates and the bound
  $2\lambda^2$ is asymptotically sharp.
\end{itemize}
For the simulations presented below,
the discretization~\eqref{eq:dCE2} has been implemented as a numerical scheme,
using the fluxes~\eqref{eq:def:diffFlux} and~\eqref{eq:def:transpFlux},
with an implicit Euler discretization in time.
\begin{figure}
	\includegraphics[width=\linewidth]{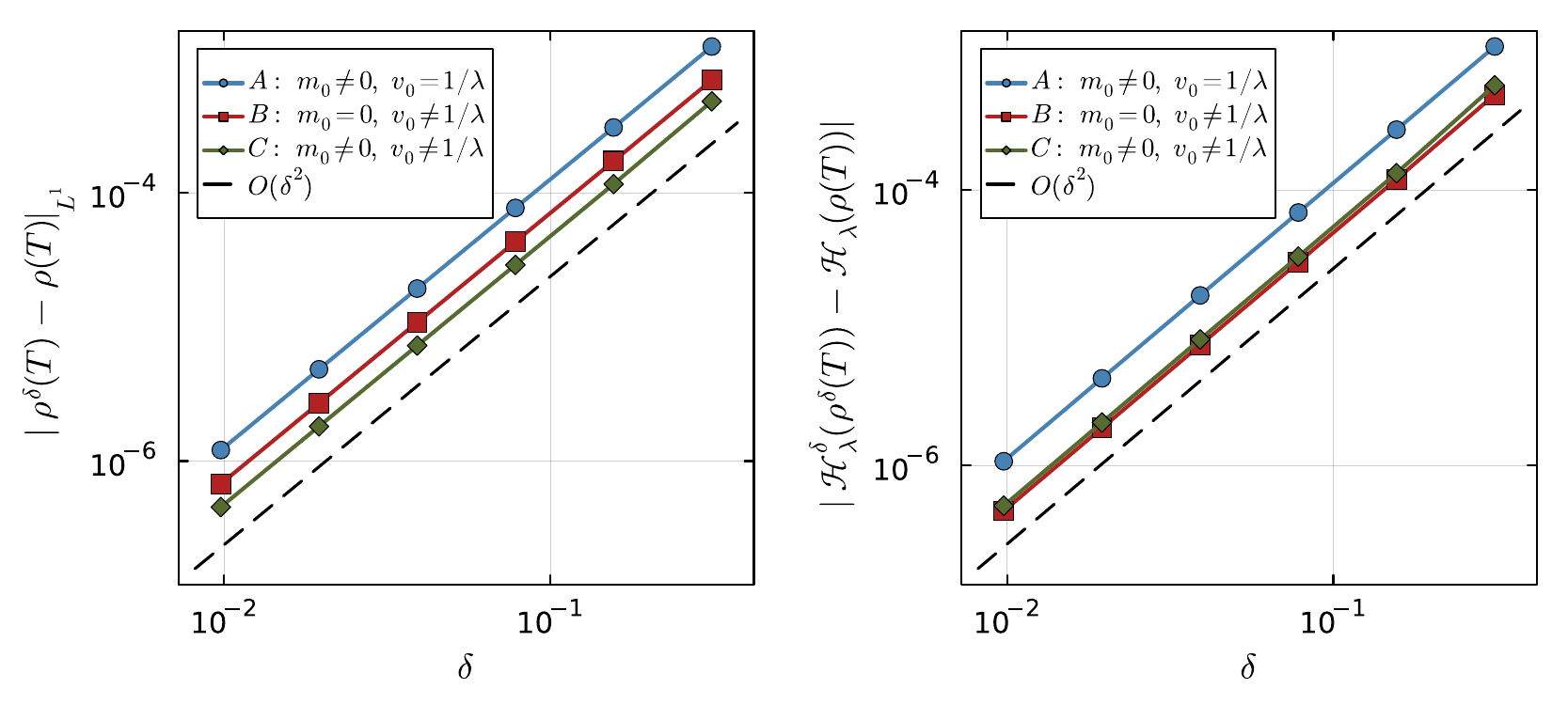}
	\caption{Numerical comparison of explicit solution from Proposition~\ref{prop:Gaussian} and solution of the scheme with parameters $\lambda=1$ and $A$: $m_0=1$, $v_0=1$; $B$: $m_0=0$, $v_0=2$; $C$: $m_0=1$, $v_0=2$. Left: the $L^1$-error; right: the difference of the exact entropy~\eqref{eq:GaussEnt} and the numerical one.}\label{fig:OrderScheme}
\end{figure}
Our method of proof for Result~\ref{res:longtime} uses the dissipation of the Otto gradient flow part alone. We expect that for centered data an improved rate of $8\lambda^2$ can be obtained by using the full dissipation, at least locally near equilibrium. The numerical results in Figure~\ref{fig:ExpCentered} confirm this expectation, showing that the entropy decays at the improved rate for centered data after a transient period. 
Notice that for off-centered data the general rate $2\lambda^2$ is attained. 
However, we note that the numerical results in Figure~\ref{fig:ExpShifted} and Figure~\ref{fig:2GaussAsym}
show non-monotone behaviour of the diffusive dissipation (see~\eqref{eq:dissipation:diff} and~\eqref{eq:dissipation:transp} in Lemma~\ref{lem:EntropyDissipation} for the exact definitions),
which is not captured by the explicit Gaussian solutions and is an interesting topic for further investigation.
We expect that for this reason the Bakry--\'{E}mery criterion is not directly generalizable to the second-order diffusive metric.

\begin{remark}[Linearization and Hermite spectral structure]\label{rem:linearization}
  The improved rate $8\lambda^2$ for centred data has a natural explanation via
  linearization of~\eqref{eq:DLSS1} around the equilibrium $\bar\rho$.
  Writing $\rho = \bar\rho(1+\varepsilon h)$ and expanding to first order in $\varepsilon$
  yields, after dividing by $\varepsilon\bar\rho$,
  \begin{align*}%
    \partial_t h = \mathcal{L}_\lambda h
    \coloneqq -\frac{1}{\bar\rho}\partial_{xx}\bigl(\bar\rho\,\partial_{xx}h\bigr)
    + \lambda\frac{1}{\bar\rho}\partial_x\bigl(\bar\rho\,\partial_x h\bigr).
  \end{align*}
  The Carrillo--Toscani splitting is inherited by the linearization
  in the form $\mathcal{L}_\lambda = -L_\lambda^2$,
  where $L_\lambda := \partial_{xx} - \lambda x\partial_x$ is the standard Ornstein-Uhlenbeck operator,
  with invariant measure $\bar\rho$.
  We refer to~\cite{AGRT} for a detailed study of such bi-Kolmogorov operators in weighted $L^2$ spaces.
  The operator~$L_\lambda$ is self-adjoint in $L^2(\bar\rho\,\dd x)$
  with Hermite polynomials as eigenfunctions,
  and eigenvalues $-n\lambda$.
  Accordingly, $\mathcal{L}_\lambda = -L_\lambda^2$ has eigenvalues $-n^2\lambda^2$.
  Explicitly:
  \begin{itemize}
  \item $n=0$: eigenvalue $0$ (mass conservation),
  \item $n=1$: eigenvalue $-\lambda^2$ (spectral gap; mean relaxation),
  \item $n=2$: eigenvalue $-4\lambda^2$ (variance relaxation),
  \item $n\geq3$: eigenvalue $-n^2\lambda^2$ (higher modes decay faster).
  \end{itemize}
  For centred perturbations, the $n=1$ mode is absent,
  so the slowest mode is $n=2$ with decay rate $-4\lambda^2$;
  the relative entropy then decays at rate $2\times 4\lambda^2 = 8\lambda^2$.
  This is consistent with the Gaussian solutions for centred data above and with the numerics in Figure~\ref{fig:ExpCentered}.
  The spectral picture also shows why the general rate $2\lambda^2$ is sharp for off-centred data:
  any $n=1$ component decays only at rate $\lambda^2$,
  and the entropy is quadratic in the perturbation.
\end{remark}

\section{Discretization}
\label{sct:discrete}
\subsection{Discretized functions and difference operator}
For a given mesh width $\delta>0$, we consider the $\delta$-spaced doubly infinite lattice~$\delta\Z\subset\R$,
i.e., elements of $\delta\Z$ are points of the form $x_k=k\delta$ with $k\in\Z$.
The dual mesh $\delta\Z^*$ consists of points $x_\kappa=\kappa\delta$
for half-integer indices $\kappa\in\Z^*:=\Z+\frac12$,
and~$x_\kappa$ is identified with the interval $(x_{\kappa-1/2},x_{\kappa+1/2})$
in between consecutive points of~$\delta\Z$.
Throughout, Latin indices such as $k$ label points of the primal grid $\delta\Z$.
Greek indices such as $\kappa$ label points of the dual grid $\delta\Z^*$.
This convention is the same as in~\cite{MRSS}.
A function $\psi:\delta\Z^*\to\R$ on the dual grid should be thought of as piecewise constant,
attaining the value $\psi_\kappa$ on $(x_{\kappa-1/2},x_{\kappa+1/2})$.
The difference quotient $\dff$ turns a function $\psi:\delta\Z^*\to\R$ on the dual grid
to a function $\dff\psi:\delta\Z\to\R$ on the primary grid,
\begin{align*}
  \dff(\psi)_k =\ \dff_k\psi = \frac{\psi_{k+1/2}-\psi_{k-1/2}}{\delta}\qquad\text{for all } k\in\Z.
\end{align*}
Accordingly, the dual operator $\dff^*$ turns a function $f:\delta\Z\to\R$ to $\dff^*f:\delta\Z^*\to\R$,
\begin{align*}
  \dff^*(f)_\kappa = \dff^*_\kappa f = \frac{f_{\kappa+1/2}-f_{\kappa-1/2}}{\delta}\qquad\text{for all } \kappa\in\Z^*,
\end{align*}
so that summation by parts reads $\delta\sum_\kappa\psi_\kappa\,\dff^*_\kappa f=-\delta\sum_k f_k\,\dff_k\psi$.
The discrete Laplace operator is defined as $\Delta:=\dff^*\dff:\delta\Z^*\to \R$, so that
\begin{align*}
  \Delta(\psi)_\kappa =\ \Delta_\kappa\psi = \frac{\psi_{\kappa+1}+\psi_{\kappa-1}-2\psi_\kappa}{\delta^2}\qquad\text{for all } \kappa\in\Z+\tfrac12.
\end{align*}
We denote by
\begin{align*}
  \dprb = \biggl\{ \rho:\delta\Z^*\to\Rnn\ \bigg|\ \delta\sum_\kappa\rho_\kappa = 1,\,\delta\sum_\kappa x_\kappa^2\rho_\kappa <\infty \biggr\}
\end{align*}
the space of all probability densities on $\delta\Z^*$ that have finite second moment.
We introduce further the discrete analogue $\bar\rho^\delta$ of the Gaussian~\eqref{eq:Gauss},
\begin{align}
  \label{eq:dGauss}
  \bar\rho^\delta_\kappa := \frac1{Z_\lambda^\delta}\exp\biggl(-\frac\lambda2x_\kappa^2\biggr)
  \quad\text{with}\quad
  Z_\lambda^\delta := \delta\sum_\kappa \exp\biggl(-\frac\lambda2x_\kappa^2\biggr),
\end{align}
i.e., $Z_\lambda^\delta>0$ is chosen such that $\bar\rho^\delta\in\dprb$.

\subsection{Discretized entropy functional}
The entropy of $\rho\in\dprb$ relative to $\bar\rho^\delta$ is given by
\begin{align}
  \label{eq:dH}
  \ent_\lambda^\delta(\rho) = \delta\sum_\kappa \biggl( \rho_\kappa\log\frac{\rho_\kappa}{\bar\rho^\delta_\kappa} - \rho_\kappa + \bar\rho^\delta_\kappa \biggr).
\end{align}
\begin{lemma}
  \label{lem:entbelow}
  For each $\lambda>0$, the relative entropy functional $\ent_\lambda^\delta$
  is minimized by~$\bar\rho^\delta$ from~\eqref{eq:dGauss}.
  More precisely, it holds
  \begin{align*}
    \ent_\lambda^\delta(\rho)\ge 0
  \end{align*}
  for each $\rho\in\dprb$, with equality iff $\rho=\bar\rho^\delta$.
\end{lemma}
\begin{proof}
  This is a classical argument:
  write the relative entropy in the form
  \begin{align}
    \label{eq:zlogz}
    \ent_\lambda^\delta(\rho)
    = \delta\sum_\kappa h\biggl(\frac{\rho_\kappa}{\bar\rho^\delta_\kappa}\biggr)\bar\rho^\delta_\kappa
    \quad\text{with}\quad
    h(z)=z\log z - z + 1.
    \end{align}
  Since $h$ is strictly convex, and since $\bar\rho^\delta$ is a probability density,
  Jensen's inequality implies that
  \begin{align*}
    \ent_\lambda^\delta(\rho)
    \ge h\biggl(\delta\sum_\kappa\frac{\rho_\kappa}{\bar\rho^\delta_\kappa}\bar\rho^\delta_\kappa\biggr)
    = h\biggl(\delta\sum_\kappa\rho_\kappa\biggr)
    = h(1) = 0,
  \end{align*}
  with equality iff $\rho_\kappa/\bar\rho^\delta_\kappa=1$ for all $\kappa$,
  i.e., iff $\rho=\bar\rho^\delta$.
\end{proof}
\begin{lemma}
  \label{lem:mom}
  For each $\lambda>0$, and for all $\rho\in\dprb$:
  \begin{align}
    \label{eq:mom}
    \frac\lambda4\delta\sum_\kappa x_\kappa^2\rho_\kappa
    \le \ent_\lambda^\delta(\rho) + \log \frac{Z_{\lambda/2}^\delta}{Z_\lambda^\delta},
  \end{align}
  with $Z^\delta$ given in~\eqref{eq:dGauss} above.
\end{lemma}
\begin{proof}
  Rewrite the entropy as follows,
  \begin{align*}
    \ent_\lambda^\delta(\rho)
    &= \delta\sum_\kappa \rho_\kappa\bigl(\log\rho_\kappa-\log\bar\rho^\delta_\kappa\bigr) \\
    &= \delta\sum_\kappa \rho_\kappa\log\rho_\kappa
      + \frac\lambda2\delta\sum_\kappa x_\kappa^2\rho_\kappa
      + \log Z_\lambda^\delta.
 \end{align*}
 And so,
 \begin{align*}
   \ent_\lambda^\delta(\rho)
   &= \frac\lambda4\delta\sum_\kappa x_\kappa^2\rho_\kappa
     + \delta\sum_\kappa\left(\rho_\kappa\log\rho_\kappa + \frac\lambda4x_\kappa^2\rho_\kappa\right)
     + \log Z_\lambda^\delta \\
   &= \frac\lambda4\delta\sum_\kappa x_\kappa^2\rho_\kappa + \ent_{\lambda/2}^\delta(\rho) + \log Z_\lambda^\delta - \log Z_{\lambda/2}^\delta.
  \end{align*}
  Using that $\ent_{\lambda/2}^\delta(\rho)\ge0$ by Lemma~\ref{lem:entbelow} above,
  we obtain~\eqref{eq:mom}.
\end{proof}
\begin{lemma}\label{lem:finiteEntropy}
  Each $\rho\in\dprb$ has finite relative entropy $\ent_\lambda^\delta(\rho)$.
\end{lemma}
\begin{proof}
  Since $\ent_\lambda^\delta$ is non-negative by Lemma~\ref{lem:entbelow} above,
  we only need to show that $\ent_\lambda^\delta(\rho)<\infty$.
  Observe that from $\rho_\kappa\ge0$ and $\delta\sum_\kappa\rho_\kappa=1$,
  it follows that $\rho_\kappa\le1/\delta$ for all $\kappa\in\delta\Z^*$.
  And so
  \begin{align*}
    \ent_\lambda^\delta(\rho)
    &= \delta\sum_\kappa \rho_\kappa\log\rho_\kappa
    - \delta\sum_\kappa \rho_\kappa\log\bar\rho^\delta_\kappa \\
    &\le \delta \sum_\kappa (-\log\delta)\rho_\kappa
    + \delta\sum_\kappa \frac\lambda2 x_\kappa^2\rho_\kappa
    + \delta\sum_\kappa (\log Z_\lambda^\delta)\rho_\kappa,
  \end{align*}
  which is finite since $\rho\in\dprb$ has finite second moment by definition.
\end{proof}
\begin{lemma}
  \label{lem:Hlsc}
  Assume that a sequence $\rho^n\in\dprb$ possesses a pointwise limit $\hat\rho\in\dprb$,
  i.e., $\rho^n_\kappa\to\hat\rho_\kappa$ for each $\kappa\in\Z^*$.
  Then
  \begin{align*}
    \ent_\lambda^\delta(\hat\rho)\le \operatorname*{liminf}_{n\to\infty}\ent_\lambda^\delta(\rho^n).
  \end{align*}
\end{lemma}
\begin{proof}
  This is an immediate consequence from Fatou's lemma:
  rewrite $\ent_\lambda^\delta$ in the form~\eqref{eq:zlogz} with the
  non-negative ``integrand'' $\kappa\mapsto h(\rho_\kappa/\bar\rho^\delta_\kappa)\bar\rho^\delta_\kappa$,
  then use the pointwise convergence.
\end{proof}

\subsection{A formal argument for consistency}
\label{sct:consistency}
We provide a short formal argument for consistency
of the discrete fluxes~$I^\delta$ and~$J^\delta$ defined in~\eqref{eq:def:diffFlux}\&\eqref{eq:def:transpFlux}
with their continuous counterparts $I$ and $J$ from~\eqref{eq:DLSS-IJ}.
Below, assume that some smooth positive density $\tilde\rho:\R\to\Rp$ is given,
and that the discrete density $\rho:\delta\Z^*\to\Rp$ is obtained by sampling,
i.e., $\rho_\kappa = \tilde\rho(x_\kappa)$ for all $\kappa\in\delta\Z^*$.

Concerning $I^\delta$:
fix some $\kappa\in\delta\Z^*$.
Then

\begin{align*}  
    I^\delta_\kappa
  &= -\frac2{\delta^2}\Bigl(e^{\frac\lambda2{\delta^2}}\sqrt{\rho_{\kappa+1}\cdot \rho_{\kappa-1}}-\rho_\kappa\Bigr) \\
  &= -\rho_\kappa \frac{2}{\delta^2} \biggl[ \exp\biggl( \log\Bigl( e^{\frac\lambda2{\delta^2}}\sqrt{\frac{\rho_{\kappa+1}\cdot \rho_{\kappa-1}}{\rho_\kappa^2}} \Bigr)\biggr) - 1 \biggr] \\
  &= -\rho_\kappa \frac{2}{\delta^2} \biggl[ \exp\biggl( \frac{\delta^2}{2} \Delta_\kappa \log \rho + \lambda \frac{\delta^2}2 \biggr) - 1 \biggr] \\
  &= -\rho_\kappa \frac{2}{\delta^2} \biggl[ \exp\biggl( \frac{\delta^2}{2} \Delta_\kappa \log \frac{\rho}{\bar\rho^\delta} \biggr) - 1 \biggr] \\
  &= -\rho_\kappa \Delta_\kappa \log \frac{\rho}{\bar\rho^\delta} + O(\delta^2) \\
  &= -\tilde\rho(x_\kappa) \,\partial_{xx} \left(\log\frac{\tilde\rho}{\bar\rho}\right)(x_\kappa) + O(\delta^2) \\
  &= I(x_\kappa) + O(\delta^2) \,.
\end{align*}
Concerning $J^\delta$:
fix some $k\in\delta\Z$.
Then
\begin{align*}
  J^\delta_{k} &= -\frac1\delta\Bigl(e^{\frac\lambda 2 \delta x_k}\rho_{k+1/2}-e^{-\frac\lambda 2 \delta x_k}\rho_{k-1/2}\Bigr) \\
  &= -\sqrt{\bar\rho^\delta_{k-1/2}\cdot \bar\rho^\delta_{k+1/2}} \ \frac1\delta\biggl( \frac{\rho_{k+1/2}}{\bar\rho^\delta_{k+1/2}} - \frac{\rho_{k-1/2}}{\bar\rho^\delta_{k-1/2}} \biggr) \\
  &= -\bar\rho(x_{k}) \partial_x \frac{\tilde\rho}{\bar\rho}(x_{k}) + O(\delta^2)
  = J(x_{k}) + O(\delta^2) \,,
\end{align*}
where we used that the central difference approximation is of order $2$.

\section{Result A: gradient flow structure}
\subsection{Gradient flow structure of DLSS}
\label{sct:cartos}
In~\cite{GST}, it has been detailed by \GS\ and collaborators
how the confined DLSS equation~\eqref{eq:DLSS1} can be understood as the metric gradient flow
of a relative Fisher information with respect to the Otto--Wasserstein distance $\wass$.
The flux representation~\eqref{eq:DLSSflux} from~\cite{CarTos}, on the other hand,
reveals an alternative formal gradient flow structure,
namely with the relative entropy $\ent_\lambda$ as potential,
and with respect to a metric $\metric$ that combines $\wass$
with the diffusive transport metric, where the latter has been recently defined in~\cite{MRSS} on the circle and in~\cite{MatthesRottSchlichting2025} on the line.
We do not give a rigorous definition of $\metric$, which is not needed in the following.
  Formally, $\metric$ is the distance induced via the Onsager (``cotangent-to-tangent'') operator $\mathcal{K}$ given in \eqref{eq:Onsager} below:
  the set of probability densities is considered as a Riemannian manifold,
  and the inverse of $\mathcal{K}$ defines a metric tensor on each tangent space.
  Distances between densities are then defined in the differential-geometric spirit as minimal length of connecting curves.
  This is analogous to the definition of $\wass$ via the Benamou--Brenier formula.
  Below, we never refer to $\metric$ explicitly, but only use the operator $\mathcal{K}$.

More specifically, introduce the symmetric and positive semi-definite (Onsager) operator~$\mathcal{K}$
at a given probability density $\rho$ by
\begin{align}
  \label{eq:Onsager}
  \mathcal{K}(\rho)\xi = \partial_{xx}\bigl(\rho\,\partial_{xx}\xi\bigr) - \lambda \partial_{x}\bigl(\rho\,\partial_{x}\xi\bigr).
\end{align}
Observing that the relative entropy $\ent_\lambda$ has variational derivative
\begin{align*}
  \ent_\lambda'(\rho) = \log\frac{\rho}{\bar\rho},
\end{align*}
it is easily deduced from the form~\eqref{eq:DLSS-IJ}  of the fluxes $I$ and $J$ that
\begin{align}
  \label{eq:DLSSgf}
  \partial_t\rho = -\mathcal{K}(\rho)\bigl[\ent_\lambda'(\rho)\bigr].
\end{align}
Combinations of two or even more different metric structures in one by addition of the Onsager operators
have been considered in a variety of cases.
Probably the most famous such compound is 
the Hellinger--Kantorovich metric $\hk$,
which combines the Otto--Wasserstein metric $\wass$ with the Hellinger distance $\hell$.
For a profound analysis of $\hk$, we refer to several works by \GS\ with Liero and Mielke~\cite{LMS1,LMS2,LMS3}.
A similar analysis for the combination $\metric$ of diffusive transport with the optimal transport metric~$\wass$
would be desirable, but is not available yet.

It does not seem that the confined DLSS equation is a contractive flow in $\metric$,
or in other words, $\ent_\lambda$ is apparently not geodesically convex in $\metric$.
The fact that the long-time behaviour of~\eqref{eq:DLSS1} can be estimated
is related to the geodesic $\lambda$-convexity of  $\ent_\lambda$ in $\wass$:
\begin{align*}
  -\frac{\dd}{\dn t}\ent_\lambda(\rho )
  &= \int\ent_\lambda'(\rho)\,\mathcal{K}(\rho)\bigl[\ent_\lambda'(\rho)\bigr]\dd x \\
  &= \int \rho\bigl[\partial_{xx}\ent_\lambda'(\rho)\bigr]^2\dd x
    + \lambda\int \rho\bigl[\partial_{x}\ent_\lambda'(\rho)\bigr]^2\dd x .
\end{align*}
The first of the integrals in the last line above is clearly non-negative;
the second one is $|\partial_{\wass}\ent_\lambda|^2(\rho)$,
the squared metric slope of the relative entropy in the $L^2$-Wasserstein distance.
By $\lambda$-convexity\footnote{The stated functional inequality is actually equivalent to Gross' logarithmic Sobolev inequality,
  which can be proven in manifold ways without any reference to $\lambda$-convexity.},
one has $|\partial_{\wass}\ent_\lambda|^2(\rho)\ge2\lambda\ent_\lambda(\rho)$, see e.g.~\cite{AGS}.
And thus one may now conclude that
\begin{align}
  \label{eq:gronwall}
  -\frac{\dd}{\dn t}\ent_\lambda(\rho )
  \ge 0 + \lambda\bigl(2\lambda\ent_\lambda(\rho)\bigr) = 2\lambda^2\ent_\lambda(\rho),
\end{align}
which yields the decay~\eqref{eq:Hdecay} by a Gronwall argument.

\subsection{Gradient flow structure for the discretization}
Our discretization~\eqref{eq:dCE2} has been chosen such that it essentially replicates the gradient flow structure above.
The main challenge in the definition of the discretization has been
to identify appropriate discrete counterparts of the operators $\partial_{xx}(\rho\,\partial_{xx}\xi)$ and $-\partial_x(\rho\,\partial_x\xi)$
inside the Onsager operator~\eqref{eq:Onsager}
such that the discrete entropy $\ent_\lambda^\delta$ is a Lyapunov function and the transport part is geodesically $\lambda^\delta$-convex
in the induced discrete Otto--Wasserstein metric, with a modulus of convexity $\lambda^\delta$ close to $\lambda$.
To make the gradient flow structure~\eqref{eq:dDLSSgf}--\eqref{eq:dOnsager} explicit,
we first express the fluxes $I^\delta$ and $J^\delta$ in terms of the relative density $u=\rho/\bar\rho^\delta$,
and then introduce mobility functions that bring the result into gradient flow form.

For a concise definition of the structure, 
recall the definition of the \emph{logarithmic mean},
\begin{equation*}
  \logmean(a,b) := \begin{cases}
    \frac{a-b}{\log a-\log b}, & a\ne b ;\\
    a , & a=b ,
  \end{cases}	 \qquad\text{for $a,b>0$,}
\end{equation*}
and introduce the \emph{geometric mean} $\geom_k(\rho):\delta\Z\to\Rnn$ of a discrete measure $\rho\in \dprb$ and the \emph{geometric Laplacian} $\geoLap_\kappa(u):\delta\Z^*\to\Rnn$ for a function $u:\delta\Z^*\to\Rnn$ by
\begin{align*}
  \geom_k(\rho) &:= \sqrt{\rho_{k+1/2}\cdot \rho_{k-1/2}}, 
  \qquad\text{ and }\qquad 
  \geoLap_\kappa(u) := \frac{2}{\delta^2}(\sqrt{u_{\kappa+1}\cdot u_{\kappa-1}}-u_\kappa).
\end{align*}
The geometric Laplacian has originally been introduced in~\cite{MRSS}
to realize the second-order chain rule~\eqref{eq:chain:Mob} below,
in parallel to the respective chain rule~\eqref{eq:chain:mob} satisfied by the logarithmic mean.
\begin{lemma}[Relative density representation]\label{lem:FluxRelDensity}
  With the relative density $u=\rho/\bar\rho^\delta$,
  the fluxes $I^\delta$~\eqref{eq:def:diffFlux} and $J^\delta$~\eqref{eq:def:transpFlux} satisfy
  \begin{align}\label{eq:FluxRelDensity}
    I^\delta_\kappa[\rho] = -\bar\rho^\delta_\kappa\,\geoLap_\kappa(u),
    \quad\text{and}\quad
    J^\delta_k[\rho] = -\geom_k(\bar\rho^\delta)\,\dff_k u.
  \end{align}
\end{lemma}
\begin{proof}
  The key ingredients are the following identities for the discrete Gaussian~\eqref{eq:dGauss},
  which follow directly from its explicit form:
  \begin{align}\label{eq:dGauss:key}
    e^{\frac{\lambda}{2}\delta^2}\sqrt{\bar\rho^\delta_{\kappa-1}\cdot \bar\rho^\delta_{\kappa+1}} = \bar\rho^\delta_\kappa,
    \qquad
    e^{\frac{\lambda}{2}\delta x_k}\bar\rho^\delta_{k+1/2}
    = e^{-\frac{\lambda}{2}\delta x_k}\bar\rho^\delta_{k-1/2}
    = \geom_k(\bar\rho^\delta).
  \end{align}
  Substituting $\rho=\bar\rho^\delta u$ into~\eqref{eq:def:diffFlux} and applying the first identity in~\eqref{eq:dGauss:key}, we arrive at
  \begin{align*}
    I_\kappa^\delta[\rho]
    &= -\frac{2}{\delta^2}\bar\rho^\delta_\kappa\Biggl(
      e^{\frac{\lambda}{2}\delta^2}
      \sqrt{\frac{\bar\rho^\delta_{\kappa-1}\cdot \bar\rho^\delta_{\kappa+1}}{(\bar\rho^\delta_\kappa)^2}}
      \sqrt{u_{\kappa-1}\cdot u_{\kappa+1}} - u_\kappa\Biggr)
    = -\bar\rho^\delta_\kappa\,\geoLap_\kappa(u)\,.
  \end{align*}
  Substituting $\rho=\bar\rho^\delta u$ into~\eqref{eq:def:transpFlux} and applying the second identity in~\eqref{eq:dGauss:key}, we obtain
  \begin{equation*}\begin{split}
    J_k^\delta[\rho]
    &= -\frac{1}{\delta}\bra*{e^{\frac{\lambda}{2}\delta x_k}\bar\rho^\delta_{k+1/2}u_{k+1/2}
      - e^{-\frac{\lambda}{2}\delta x_k}\bar\rho^\delta_{k-1/2}u_{k-1/2}} \\
    &= -\frac{\geom_k(\bar\rho^\delta)}{\delta}\bra*{u_{k+1/2} - u_{k-1/2}}
    = -\geom_k(\bar\rho^\delta)\,\dff_k u\,,
  \end{split}\end{equation*}
  which concludes the proof of~\eqref{eq:FluxRelDensity}.
\end{proof}
From the identities~\eqref{eq:FluxRelDensity} we obtain our choice of mobility functions.
By using the chain rules
\begin{align}\label{eq:chain:Mob}
  \geoLap_\kappa(u) &= \logmean\bra*{\sqrt{u_{\kappa+1}\cdot u_{\kappa-1}},u_\kappa}\,\Delta_\kappa\log u \\
\intertext{and}
\label{eq:chain:mob}
  \dff_k u &= \logmean(u_{k+1/2},u_{k-1/2})\,\dff_k\log u\,,
\end{align}
the right-hand sides in~\eqref{eq:FluxRelDensity} take gradient flow form once we define
mobility functions $\mob^\delta$ and $\Mob^\delta$ on the primal and dual grid, in terms of $u=\rho/\bar\rho^\delta$, as follows
\begin{align}\label{eq:def:Mob}
\Mob^\delta_\kappa(\rho) := \bar\rho^\delta_\kappa\logmean(\sqrt{u_{\kappa+1}\cdot u_{\kappa-1}},u_\kappa)
= \logmean\bigl(e^{\tfrac{\lambda}2 \delta^2}\sqrt{\rho_{\kappa+1}\cdot \rho_{\kappa-1}},\rho_\kappa\bigr) \,,
\end{align}
where we used the identity~\eqref{eq:dGauss:key} for the discrete Gaussian.
Likewise, by a similar calculation, we have
\begin{align}\label{eq:def:mob}
  \mob^\delta_k(\rho) := \geom_k(\bar\rho^\delta)\logmean(u_{k+1/2},u_{k-1/2}) = \logmean\Bigl(e^{\tfrac\lambda2 \delta x_{k}}\rho_{k+1/2}, e^{-\tfrac\lambda2 \delta x_{k}}\rho_{k-1/2}\Bigr) \,.
\end{align}
Hence, the chain rules~\eqref{eq:chain:Mob}--\eqref{eq:chain:mob} immediately yield
the following gradient flow representation~\eqref{eq:dDLSSgf} of the flux formulation~\eqref{eq:dCE2} of the confined DLSS equation, which concludes Result~\ref{res:GF}.
\begin{lemma}[Gradient flow representation]\label{lem:FluxMob}
  The fluxes $I^\delta$~\eqref{eq:def:diffFlux} and $J^\delta$~\eqref{eq:def:transpFlux},
  with mobilities $\Mob^\delta$~\eqref{eq:def:Mob} and $\mob^\delta$~\eqref{eq:def:mob}, satisfy
  \begin{align}
    \label{eq:FluxMob}
    I^\delta_\kappa[\rho] = -\Mob^\delta_\kappa(\rho)\,\Delta_\kappa\bigl(\ent_\lambda^\delta\bigr)'(\rho),
    \quad\text{and}\quad
    J^\delta_k[\rho] = -\mob^\delta_k(\rho)\,\dff_k\bigl(\ent_\lambda^\delta\bigr)'(\rho).
  \end{align}
\end{lemma}
\begin{proof}
  Since $(\ent_\lambda^\delta)'(\rho)_\kappa = \log u_\kappa$,
  the identities follow from Lemma~\ref{lem:FluxRelDensity}
  and the chain rules~\eqref{eq:chain:Mob}--\eqref{eq:chain:mob}.
\end{proof}
On a purely formal level, the long-time behaviour of~\eqref{eq:dDLSSgf} can now be obtained
in complete analogy to the continuous case~\eqref{eq:DLSSgf}.
Summing by parts, one finds
\begin{align*}
  -\frac{\dd}{\dn t}\ent_\lambda^\delta(\rho )
  &= \delta\sum_\kappa\bigl[(\ent_\lambda^\delta)'(\rho)\bigr]_\kappa \,\mathcal{K}^\delta_\kappa(\rho)\bigl[(\ent_\lambda^\delta)'(\rho)\bigr] \\
  &= \delta\sum_\kappa \Mob^\delta_\kappa(\rho)\bigl[\Delta_\kappa\bigl((\ent_\lambda^\delta)'(\rho)\bigr)\bigr]^2
    + \lambda\delta\sum_k \mob^\delta_k(\rho)\bigl[\dff_k\bigl((\ent_\lambda^\delta)'(\rho)\bigr)\bigr]^2 .
 \end{align*}
The first sum in the last line above is non-negative,
the second sum represents the squared slope of the relative entropy $\ent_\lambda^\delta$
in the $\mob^\delta$-induced metric.
And the results from~\cite{ErbarMaas2012,FathiMaas2016,Mielke2013} for \emph{finite} Markov chains
are expected to imply that the relative entropy $\ent_\lambda^\delta$ is $\lambda^\delta$-convex
with respect to the generalized transport metric induced via the mobility $\mob^\delta$,
and with some $\lambda^\delta\approx\lambda$.
This would allow us to proceed to the discrete analogue of~\eqref{eq:gronwall}.

For the rigorous argument presented in this paper,
we did not try to extend the aforementioned convexity theory
from finite to countably infinite Markov chains.
Instead, for the passage to the Gronwall estimate,
we use a modified logarithmic Sobolev inequality from~\cite{CaputoPraPosta2009} in Section~\ref{sec:mLSI} below.

\section{Result B: global existence and positivity}
Let $\delta>0$ be fixed.
Our strategy is to obtain the time-global non-negative solution to the initial value problem on $\delta\Z$
as a limit of strictly positive solutions on finite grids.
The finite size allows us to apply standard existence results for ODEs,
and the strict positivity leads to smoothness of these approximations.

\subsection{Finite grid and boundary conditions}
We approximate solutions to~\eqref{eq:scheme} on the finite grid $\delta\llbracket-K,K\rrbracket$ with $2K$ cells
by restricting the continuity equation~\eqref{eq:dCE2} to this domain.
The discrete Laplacian $\Delta$ in the diffusive term of~\eqref{eq:dCE2} is interpreted as the
reflected Laplacian on $\delta\llbracket-K,K\rrbracket$, and the discrete divergence $\dff^*$
in the transport term is subjected to a no-flux boundary condition. 
The choice is motivated from the structure of the Onsager operator $\mathcal{K}^\delta$ in~\eqref{eq:dOnsager}, which is a sum of an Otto--Wasserstein part and a diffusive transport part, for which the former is naturally associated with no-flux boundary conditions, whereas the latter is naturally associated with reflected Laplacian boundary conditions.
Together these boundary conditions on the fluxes read
\begin{align}
  \label{eq:newbc}
  \dff_{K} I^\delta[\rho] = \dff_{- K} I^\delta[\rho]=0
  \qquad \text{and} \qquad
  J^\delta_K[\rho] = J^\delta_{-K}[\rho] = 0.
\end{align}
The no-flux condition $J_{\pm K}[\rho]=0$ determines the values of $\rho$ at the ghost points $\kappa=\pm(K+1/2)$,
\begin{align}
  \label{eq:ghost}
  \rho_{K+1/2} = \rho_{K-1/2} e^{-\lambda \delta x_K}
  \quad\text{and}\quad
  \rho_{-K-1/2} = \rho_{-K+1/2} e^{-\lambda \delta x_K},
\end{align}
which allows us to express the diffusive flux near the boundary explicitly as
\begin{align}
  \label{eq:bcm:diffFlux}
  I^\delta_{-K+1/2}[\rho]
  &=- \frac{2}{\delta^2}\sqrt{\rho_{-K+1/2}} \bigl( e^{\frac{\lambda}{2} \delta x_{-K+1}} \sqrt{\rho_{-K+3/2}} - \sqrt{\rho_{-K+1/2}}\bigr) \,, \\
  \label{eq:bcp:diffFlux}
  I^\delta_{K-1/2}[\rho]
  &=- \frac{2}{\delta^2}\sqrt{\rho_{K-1/2}} \bigl( e^{-\frac{\lambda}{2} \delta x_{K-1}} \sqrt{\rho_{K-3/2}} - \sqrt{\rho_{K-1/2}}\bigr) \,. 
\end{align}
The reflected Laplacian condition in~\eqref{eq:newbc} further yields ghost-point flux values
$I^\delta_{K+1/2}[\rho]=I^\delta_{K-1/2}[\rho]$ and $I^\delta_{-K-1/2}[\rho]=I^\delta_{-K+1/2}[\rho]$.

\subsection{Local existence on a finite grid}
\begin{lemma}
  \label{lem:sol1}
  Given a positive initial datum $\rho^0\in\dprbK$,
  there exists a unique maximal positive solution $\rho:[0,T)\to\dprbK$
  to~\eqref{eq:scheme} subject to~\eqref{eq:newbc} on an interval $[0,T)$.
  That solution is uniformly bounded,
  \begin{align}
    \label{eq:inftybound}
    \rho_\kappa(t) \le \delta^{-1}.
  \end{align}
  And either $T=+\infty$, or positivity is lost at $T$ in the sense that
  \begin{align}
    \label{eq:breakdown}
    \liminf_{t\nearrow T}\min_\kappa \rho_\kappa(t)=0.
  \end{align}
\end{lemma}
\begin{proof}
  The Cauchy--Lipschitz theorem applies since the right-hand side of~\eqref{eq:dCE2} defines
  a continuously differentiable (thus locally Lipschitz) vector field
  on the open set of strictly positive densities $\rho\in\R^{2K}$.
  Thus, there is a maximal solution $\rho:[0,T)\to\R^{2K}$,
  and $T=+\infty$ unless $\rho(t)$ ``reaches the boundary'' as $t\nearrow T$,
  i.e.,
  \[
    \limsup_{t\nearrow T}\max_\kappa\rho_\kappa(t)=\infty
    \quad\text{or}\quad
    \liminf_{t\nearrow T}\min_\kappa \rho_\kappa(t)=0.   
  \]
  To verify that $\rho(t)\in\dprbK$ for each $t\in(0,T)$, it suffices to show that mass is conserved:
  the time derivative of the total mass along~\eqref{eq:dCE2}
  vanishes thanks to the boundary conditions~\eqref{eq:newbc}:
  \begin{align*}
    \frac{\dd}{\dn t}&\,\delta\sum_{|\kappa|< K} \rho_\kappa
    = \delta\sum_{|\kappa|<K}\Delta_\kappa I[\rho]
      - \lambda\delta\sum_{|\kappa|<K}\dff^*_\kappa J[\rho] \\
    &= \sum_{|\kappa|< K}\bigl(\dff_{\kappa+1/2}I[\rho]-\dff_{\kappa-1/2}I[\rho]\bigr)
      - \lambda\sum_{|\kappa|<K}\bigl(J_{\kappa+1/2}[\rho]-J_{\kappa-1/2}[\rho]\bigr) \\
    &= \bigl(\dff_KI[\rho]-\dff_{-K}I[\rho]\bigr) - \lambda \bigl(J_K[\rho]-J_{-K}[\rho]\bigr)
    = 0.
  \end{align*}
  Mass conservation implies that, as long as the solution remains positive, it also remains bounded:
  for every $\kappa^*$,
  \begin{align*}
    \delta \rho_{\kappa^*}(t) \le \delta\sum_\kappa \rho_\kappa(t) = \delta\sum_\kappa\rho_\kappa^0 = 1,
  \end{align*}
  which further implies~\eqref{eq:inftybound}.
  The only possibility for break-down of the solution at finite time is thus loss of positivity~\eqref{eq:breakdown} for $t\nearrow T$.
\end{proof}
\begin{lemma}
  \label{lem:tLipschitz}
  The maximal solution $\rho$ obtained in Lemma~\ref{lem:sol1} above is globally Lipschitz in time.
  More specifically, for each $\kappa^*\in\Z^*$ with $|\kappa^*|<K$, and all $t,t'\in[0,T)$,
  we have
  \begin{align}
    \label{eq:tLipschitz}
    \bigl|\rho_{\kappa^*}(t)-\rho_{\kappa^*}(t')\bigr| \le L_{\kappa^*}|t-t'|
    \quad \text{with} \quad L_{\kappa^*}:=16e^{\frac\lambda2\delta^2}\delta^{-5} + 4\lambda e^{\frac\lambda2\delta^2(|\kappa^*|+1/2)}\delta^{-3}.
  \end{align}
\end{lemma}
\begin{proof}
  First, let us assume that $|\kappa^*|<K-2$.
  Directly from the representation~\eqref{eq:dCE2},
  we shall estimate $\rho_{\kappa^*}$'s time derivative as $|\dot\rho_{\kappa^*}|\le L_{\kappa^*}$,
  with $L_{\kappa^*}$ given in~\eqref{eq:tLipschitz}.
  Recalling the absolute bound~\eqref{eq:inftybound},
  we obtain from~\eqref{eq:def:diffFlux} and~\eqref{eq:def:transpFlux}, respectively, that
  \begin{align*}
    \bigl|I_\kappa[\rho]\bigr|
    &\le \frac{2}{\delta^2}\bigl(e^{\frac\lambda2\delta^2}\sqrt{\rho_{\kappa+1}\cdot \rho_{\kappa-1}}+\rho_\kappa\bigr)
      \le  4e^{\frac\lambda2\delta^2}\delta^{-3}, \\
    \bigl|J_k[\rho]\bigr|
    &\le \frac1\delta\bigl(e^{\frac\lambda2\delta x_k}\rho_{k+1/2}+e^{-\frac\lambda2\delta x_k}\rho_{k-1/2}\bigr)
      \le 2e^{\frac\lambda2\delta|x_k|}\delta^{-2}.
  \end{align*}
  Hence, it follows from~\eqref{eq:dCE2} via the definitions of the discrete operators $\Delta$ and $\dff^*$ that
  \begin{align*}
    \bigl|\dot\rho_{\kappa^*}\!\bigr|
    &\le \bigl|\Delta_{\kappa^*} I[\rho]\bigr| + \lambda \bigl|\dff^*_{\kappa^*} J[\rho]\bigr| \\
    &\le \frac{1}{\delta^2}\Bigl(2\bigl|I_{\kappa^*}[\rho]\bigr|+\bigl|I_{\kappa^*+1}[\rho]\bigr|+\bigl|I_{\kappa^*-1}[\rho]\bigr|\Bigr)
    + \frac{\lambda}\delta \Bigl(\bigl|J_{\kappa^*+1/2}[\rho]\bigr|+\bigl|J_{\kappa^*-1/2}[\rho]\bigr|\Bigr)\\
    &\le L_{\kappa^*},
  \end{align*}
  where we have used that $|x_{\kappa^*\pm1/2}|\le\delta(|\kappa^*|+1/2)$.

  If $\kappa^*=\pm(K-3/2)$ or $\kappa^*=\pm(K-1/2)$,
  the alternative forms of $I_\kappa[\rho]$ from~\eqref{eq:bcm:diffFlux} and~\eqref{eq:bcp:diffFlux} and for~$J_k[\rho]$ from~\eqref{eq:newbc} need to be used;
  it is easily seen that the estimate~\eqref{eq:tLipschitz} also holds at those sites.
\end{proof}

\subsection{Preservation of positivity}
We are now going to show that the necessary criterion~\eqref{eq:breakdown} for finite-time blow-up of solutions is never met.
\begin{proposition}
  \label{prp:positive}
  If $\rho:[0,T)\to\dprbK$ is the solution to~\eqref{eq:scheme} subject to~\eqref{eq:newbc}
  with positive initial datum, then~\eqref{eq:breakdown} is not true.
  Consequently, the maximal solution to the initial value problem is defined for all positive times, i.e., $T=+\infty$.
\end{proposition}
\begin{proof}
  Suppose that a solution satisfies~\eqref{eq:breakdown} for a finite terminal time $T\in(0,\infty)$.
  By uniform Lipschitz continuity in time, see Lemma~\ref{lem:tLipschitz},
  there is a limiting profile $\hat\rho\in\dprbK$
  such that
  \begin{align*}
    \lim_{t\nearrow T}\rho_\kappa(t) = \hat\rho_\kappa \quad \text{at each $\kappa$}.
  \end{align*}
  The assumed break-down~\eqref{eq:breakdown} implies that there is a $\kappa^*$ such that $\hat\rho_{\kappa^*}=0$.
  On the other hand, as a probability density, $\hat\rho$ cannot be identically zero.
  So $\kappa^*$ can be chosen such that ---
  without loss of generality\footnote{If $\kappa^*=K-1/2$ is on the right boundary,
    consider the situation obtained by reflection $\kappa\leftrightarrow-\kappa$ instead.}
  --- $\hat\rho_{\kappa^*+1}>0$.

  Since $(0,T)\ni t\mapsto\rho_{\kappa^*}(t)$ is continuously differentiable and positive,
  the mean value theorem implies the existence of an increasing sequence $(t_n)\subset(0,T)$ with $\lim_{n\to\infty}t_n=T$
  such that $\partial_t\rho_{\kappa^*}(t_n)<0$ for each $n$.
  Below, we are now going to show that $\lim_{n\to\infty}\partial_t\rho_{\kappa^*}(t_n)>0$, which is a contradiction.

  Recall the representation of $\dot\rho$ in~\eqref{eq:dCE2}.
  First, assume that $x_{\kappa^*}$ is not near the boundary, i.e., $|\kappa^*|<K-2$.
  Then $I_\kappa[\rho]$ for $\kappa\in\{\kappa^*+1,\kappa^*,\kappa^*-1\}$ is given by~\eqref{eq:def:diffFlux},
  and $J_k[\rho]$ for $k\in\{\kappa^*+1/2,\kappa^*-1/2\}$ is given by~\eqref{eq:def:transpFlux}.
  By continuity of these expressions, $I_\kappa[\rho(t_n)]\to I_\kappa[\hat\rho]$ and $J_k[\rho(t_n)]\to J_k[\hat\rho]$.
  Using $\hat\rho_{\kappa^*}=0$ and $\hat\rho_{\kappa^*+1}>0$, we obtain
  \begin{align}
    \label{eq:I1}
    I^\delta_{\kappa^*}[\hat\rho] &= -\frac{2e^{\frac\lambda2\delta^2}}{\delta^2}\sqrt{\hat\rho_{\kappa^*+1}\cdot \hat\rho_{\kappa^*-1}} \le 0, \\
    \label{eq:I2}
    I^\delta_{\kappa^*+1}[\hat\rho] &= \frac2{\delta^2}\hat\rho_{\kappa^*+1} > 0, \\
    \label{eq:I3}
    I^\delta_{\kappa^*-1}[\hat\rho] &= \frac2{\delta^2}\hat\rho_{\kappa^*-1} \ge 0,     
  \end{align}
  and consequently
  \begin{align}
    \label{eq:help1}
    \Delta_{\kappa^*}I^\delta[\hat\rho] =  \frac1{\delta^2}\bigl(I^\delta_{\kappa^*+1}[\hat\rho]+ I^\delta_{\kappa^*-1}[\hat\rho]- 2I^\delta_{\kappa^*}[\hat\rho]\bigr) > 0.
  \end{align}
  Similarly,
  \begin{align}
    \label{eq:J1}
    J^\delta_{\kappa^*+1/2}[\hat\rho] &= -\frac{e^{\frac\lambda2\delta x_{\kappa^*+1/2}}}\delta \hat\rho_{\kappa^*+1} < 0, \\
    \label{eq:J2}
    J^\delta_{\kappa^*-1/2}[\hat\rho] &= \frac{e^{-\frac\lambda2\delta x_{\kappa^*-1/2}}}\delta \hat\rho_{\kappa^*-1} \ge  0, 
  \end{align}
  which implies that
  \begin{align}
    \label{eq:help2}
    -\lambda\dff_{\kappa^*}^*J^\delta[\hat\rho] 
    = \frac\lambda\delta\bigl( J^\delta_{\kappa^*-1/2}[\hat\rho] - J^\delta_{\kappa^*+1/2}[\hat\rho] \bigr) > 0.
  \end{align}
  In summary, we obtain that
  \begin{align}
    \label{eq:pcontradiction}
    \lim_{n\to\infty}\partial_t\rho_{\kappa^*}(t_n)
    = \Delta_{\kappa^*}I^\delta[\hat\rho] - \lambda\dff_{\kappa^*}^*J^\delta[\hat\rho]
    > 0.
  \end{align}
  This is the desired contradiction in this case.

  Next, suppose that $\kappa^*=-K+3/2$.
  Then $J^\delta_{\kappa^*\pm1/2}[\hat\rho]$ are as in~\eqref{eq:J1}\&\eqref{eq:J2} above, leading to~\eqref{eq:help2}.
  Also $I^\delta_{\kappa^*}[\hat\rho]$ and $I_{\kappa^*+1}[\hat\rho]$ are as in~\eqref{eq:I1}\&\eqref{eq:I2};
  the only change is that $I_{\kappa^*-1}[\hat\rho]$ now needs to be evaluated according to~\eqref{eq:bcm:diffFlux},
  \begin{align*}
    I^\delta_{\kappa^*-1}[\hat\rho] = \frac2{\delta^2}\hat\rho_{-K+1/2}\ge 0,
  \end{align*}
  which is consistent with the ``$\ge0$'' obtained in~\eqref{eq:I3} and
  thus leads to the same contradiction~\eqref{eq:pcontradiction}.

  If instead $\kappa^*=K-3/2$,
  then $J^\delta_{\kappa^*\pm1/2}[\hat\rho]$ are still unchanged, and so is~\eqref{eq:help2}.
  $I^\delta_{\kappa^*}[\hat\rho]$ and $I_{\kappa^*-1}[\hat\rho]$ remain the same as above,
  and instead of~\eqref{eq:I2}, one uses~\eqref{eq:bcp:diffFlux} to find
  \begin{align*}
    I^\delta_{\kappa^*+1}[\hat\rho] = \frac2{\delta^2}\hat\rho_{K-1/2} > 0,    
  \end{align*}
  which still implies~\eqref{eq:help1},
  so~\eqref{eq:pcontradiction} holds also in this case.
  
  Finally, suppose $\kappa^*=-K+1/2$ is directly at the left boundary.
  Then $J^\delta_{\kappa^*+1/2}[\hat\rho]$ is as in~\eqref{eq:J1},
  while $J^\delta_{\kappa^*-1/2}[\hat\rho]=0$ by the boundary condition~\eqref{eq:newbc},
  which is still consistent with the ``$\ge0$'' in~\eqref{eq:J2}.
  While~\eqref{eq:I2} remains unchanged,
 ~\eqref{eq:I1} is replaced by $I^\delta_{\kappa^*}[\hat\rho] = 0$ because of~\eqref{eq:bcm:diffFlux},
  and with the boundary condition $\dff_{-K}I^\delta[\rho]=0$ from~\eqref{eq:newbc}, so that
  \begin{align*}
    \Delta_{\kappa^*}I^\delta[\hat\rho]
    = \frac1\delta\bigl(\dff_{\kappa^*+1/2}I[\hat\rho]-\dff_{\kappa^*-1/2}I[\hat\rho]\bigr)
    = \frac1{\delta^2}I_{\kappa^*+1}[\hat\rho]>0.
  \end{align*}
  Again,~\eqref{eq:pcontradiction} follows.
\end{proof}

\subsection{Entropy dissipation}
\begin{lemma}
  \label{lem:EntropyDissipation}
  Along the solutions $\rho:[0,\infty)\to\dprbK$ considered above,
  the entropy is dissipated as follows:
  \begin{align}\label{eq:EntropyDissipation}
    -\frac{\dd}{\dn t}\ent_\lambda^\delta(\rho)
    = \mathcal{D}^\delta_{\mathrm{diff}}(\rho) + \lambda \mathcal{D}^\delta_{\mathrm{transp}}(\rho) \,,
  \end{align}
  with the dissipations given by
  \begin{align}
  	\mathcal{D}^\delta_{\mathrm{diff}}(\rho) &:= \delta\sum_{|\kappa|\le K-1/2} \Mob^\delta_\kappa(\rho)\Bigl|\Delta_\kappa \bigl(\ent_\lambda^\delta\bigr)'(\rho)\Bigr|^2 \label{eq:dissipation:diff} \\
  	\mathcal{D}^\delta_{\mathrm{transp}}(\rho) &:= \delta\sum_{|k|\le K-1} \mob^\delta_k(\rho)\Bigl|\dff_k \bigl(\ent_\lambda^\delta\bigr)'(\rho)\Bigr|^2
  	\label{eq:dissipation:transp}
  \end{align}
  and with the understanding that, see~\eqref{eq:bcm:diffFlux}--\eqref{eq:bcp:diffFlux},
  \begin{align}
    \label{eq:interpretation}
    \begin{split}
      &\Mob^\delta_{\pm(K-1/2)}(\rho)\Bigl[\Delta_{\pm(K-1/2)}\bigl(\ent_\lambda^\delta\bigr)'(\rho)\Bigr]^2 \\
      &\qquad := \pm\frac1\delta\, \dff_{\pm(K-1)}\left(\bigl(\ent_\lambda^\delta\bigr)'(\rho)\right)\, I^\delta_{\pm(K-1/2)}[\rho] \ge 0.
      \end{split}
  \end{align}
  In particular, $\ent_\lambda^\delta(\rho)$ is monotonically decreasing in time.
\end{lemma}
Notice that~\eqref{eq:interpretation} is coherent with the values~\eqref{eq:ghost}
that are implied at the ``ghost points'' $\kappa=\pm(K+1/2)$ by the boundary conditions~\eqref{eq:newbc}.
\begin{proof}
  Since $\rho$ is positive, the time derivative of $\ent_\lambda^\delta(\rho)$
  can be directly obtained by differentiation of the logarithmic term in the sum,
  \begin{align*}
    \frac{\dd}{\dn t}\ent_\lambda^\delta(\rho)
    = \delta\sum_{|\kappa|\le K-1/2} \bigl(\ent_\lambda^\delta\bigr)'(\rho)_\kappa\,\dot\rho_\kappa.
  \end{align*}
  We substitute the evolution equation~\eqref{eq:dCE2} and sum by parts twice, using the
  boundary conditions $\dff_{\pm K}I^\delta[\rho]=0$ and $J^\delta_{\pm K}[\rho]=0$ from~\eqref{eq:newbc}:
  \begin{align*}
    \MoveEqLeft -\frac{\dd}{\dn t}\ent_\lambda^\delta(\rho)
    = +\delta\sum_{|k|\le K-1}\dff_k\bigl(\ent_\lambda^\delta\bigr)'(\rho)\cdot \dff_kI[\rho]
      - \lambda\delta\sum_{|k|\le K-1}\dff_k\bigl(\ent_\lambda^\delta\bigr)'(\rho)\cdot J_k[\rho] \\
    &= -\delta\sum_{|\kappa|\le K-3/2}\Delta_\kappa\bigl(\ent_\lambda^\delta\bigr)'(\rho)\cdot I_\kappa[\rho]
      + \dff_{K-1}\bigl(\ent_\lambda^\delta\bigr)'(\rho)\cdot I_{K-1/2}[\rho] \\
    &\qquad - \dff_{-(K-1)}\bigl(\ent_\lambda^\delta\bigr)'(\rho)\cdot I_{-(K-1/2)}[\rho]
      - \lambda\delta\sum_{|k|\le K-1}\dff_k\bigl(\ent_\lambda^\delta\bigr)'(\rho)\cdot J_k[\rho].
  \end{align*}
  Now replace, in accordance with Lemma~\ref{lem:FluxMob}, specifically~\eqref{eq:FluxMob},
  \begin{align*}
    I_\kappa[\rho] = -\Mob^\delta_\kappa(\rho)\,\Delta_\kappa\bigl(\ent_\lambda^\delta\bigr)'(\rho)
    \qquad \text{and}\qquad
    J_k[\rho] = -\mob^\delta_k(\rho)\,\dff_k\bigl(\ent_\lambda^\delta\bigr)'(\rho)
  \end{align*}
  for $|\kappa|\le K-3/2$ and $|k|\le K-1$,
  and apply the interpretation~\eqref{eq:interpretation};
  this yields~\eqref{eq:EntropyDissipation}.
  For proving that the products in~\eqref{eq:interpretation} are actually non-negative,
  it suffices to compare the signs of differences of the square roots in the expressions~\eqref{eq:def:diffFlux} and~\eqref{eq:def:transpFlux}
  with the sign of the difference quotients $\dff_{\pm(K-1)}\bigl(\ent_\lambda^\delta\bigr)'(\rho)$.
\end{proof}
\begin{corollary}
  \label{cor:tight}
  There is a $K$-independent constant $C$ such that the second moment
  of the solution $\rho$ is controlled as follows:
  \begin{align*}
    \frac\lambda4 \delta\sum_\kappa x_\kappa^2\rho_\kappa \le \ent_\lambda^\delta(\rho^0) + C.
  \end{align*}
\end{corollary}
\begin{proof}
  This follows by combining the time-monotonicity of $\ent_\lambda^\delta(\rho)$
  with the estimate~\eqref{eq:mom} on the second moment from Lemma~\ref{lem:mom}.
\end{proof}
\subsection{Passage to the infinite grid}
\label{sct:proofB}
We are now in the position to obtain the result on global existence and positivity.
\begin{proof}[Proof of Result B]
  Let an initial datum $\rho^0\in\dprb$ be given, then thanks to Lemma~\ref{lem:finiteEntropy} it is of finite entropy.
  For each $K=2,3,\ldots$, define a strictly positive $\rho^{K,0}\in\dprbK$
  such that $\rho^{K,0}_\kappa\to\rho^0_\kappa$ for each $\kappa$,
  and such that
  \begin{align}
    \label{eq:Hlimit}
    \lim_{K\to\infty}\ent_\lambda^\delta(\rho^{K,0})  = \ent_\lambda^\delta(\rho^0).
  \end{align}
  Let $\rho^K:[0,\infty)\to\dprbK$ be the unique and positive global solution
  to~\eqref{eq:dCE2} with initial datum $\rho^{K,0}$.

  Fix some $\kappa^*\in\Z^*$.
  The curves $t\mapsto\rho^K_{\kappa^*}(t)$ are Lipschitz continuous by~\eqref{eq:tLipschitz},
  with the constant $L_{\kappa^*}$ that is uniform with respect to $t\ge0$ and independent of $K>|\kappa^*|$.
  By the Arzelà--Ascoli theorem and a diagonal argument (with respect to the compact time interval $[0,T]$),
  any sequence of these curves contains a subsequence that is uniformly convergent on each $[0,T]$. 

  By means of another diagonal argument (with respect to the countably many grid sites $\kappa^*$),
  we thus obtain a sequence of $K_n\to\infty$
  and a corresponding limit curve $\rho^*:[0,\infty)\to\R^{\delta\Z^*}$
  such that, for each $\kappa\in\Z^*$,
  the curve $\rho^*_\kappa:[0,\infty)\to\R$ is globally Lipschitz continuous,
  and is the locally uniform limit of the subsequence $\rho^{K_n}_\kappa$.
  It then follows that also the fluxes $I^\delta_\kappa[\rho^{K_n}]$ and $J^\delta_k[\rho^{K_n}]$ from~\eqref{eq:def:diffFlux} and~\eqref{eq:def:transpFlux}
  --- each of which depends on the respective three curves $\rho^{K_n}_\kappa$, $\rho^{K_n}_{\kappa+1}$ and $\rho^{K_n}_{\kappa-1}$ through continuous operations ---
  converge locally uniformly as well.
  And since $\rho^{K_n}_\kappa$ satisfies the differential equation~\eqref{eq:dCE2},
  the curve $\rho^*_\kappa$ is continuously differentiable and still satisfies~\eqref{eq:dCE2}.

  Next, we show that $\rho^*$ is actually a curve in $\dprb$.
  Non-negativity is an immediate consequence of the positivity of the $\rho^{K_n}$.
  Concerning the normalization:
  let $t\ge 0$ be fixed and let $\hat\rho^{K_n}(t)$ be the extension of $\rho^{K_n}(t)$ to $\delta\Z^*$ by zero.
  Then $\hat\rho^{K_n}_\kappa(t)\to\rho^*_\kappa(t)$ for each $\kappa$,
  and the tightness estimate from Corollary~\ref{cor:tight} implies
  that $\hat\rho^{K_n}(t)$ has $n$-uniformly bounded second moment.
  This in combination implies that the normalization
  \begin{align*}
    \delta \sum_{\kappa}\hat\rho^{K_n}_\kappa(t) = 1
  \end{align*}
  is inherited by the limit $\rho^*(t)$. 

  It remains to be shown that $\rho^*$ becomes instantaneously positive.
  To that end, suppose that $\rho^*_{\kappa^*}(t^*)=0$ at some $t^*>0$ and $\kappa^*\in\Z^*$.
  Since $\rho^*_{\kappa^*}(t)\ge0$ for all $t\ge0$, differentiability at $t=t^*$ implies
  \begin{align}
    \label{eq:againpositive}
    \partial_t\rho^*_{\kappa^*}(t^*) = 0.
  \end{align}
  Below, we produce a contradiction by showing that $\partial_t\rho^*_{\kappa^*}(t^*)>0$.
  
  Since $\rho^*(t^*)\in\dprb$, we cannot have $\rho^*_\kappa(t^*)=0$ for all $\kappa\in\delta\Z^*$,
  and so, without loss of generality, we may assume that $\rho^*_{\kappa^*+1}(t^*)>0$.
  Now, on the grounds of $\rho^*_{\kappa^*}(t^*)=0<\rho^*_{\kappa^*+1}(t^*)$,
  we can directly adapt the argument from the proof of Proposition~\ref{prp:positive},
  with $\rho^*(t^*)$ in place of $\hat\rho$ to conclude~\eqref{eq:pcontradiction},
  i.e., that
  \begin{align*}
    \partial_t\rho^*_{\kappa^*}(t^*) = \Delta_{\kappa^*}I^\delta[\rho^*(t^*)] - \lambda \dff^*_{\kappa^*}J^\delta[\rho^*(t^*)] > 0,
  \end{align*}
  which contradicts~\eqref{eq:againpositive}.
\end{proof}

\section{Result C: long-time asymptotics via a discrete modified logarithmic Sobolev inequality}
\label{sec:mLSI}
We first prove the decay of the relative entropy for each approximate solution on a finite grid
and then pass to the infinite-grid limit.

The key result driving the exponential convergence is the following
entropy--entropy-dissipation inequality.
\begin{proposition}[Discrete modified log-Sobolev inequality]\label{prop:dLSI}
	For any $\lambda>0$, $\delta>0$, and $\rho\in\dprbK$ it holds
	\begin{equation}\label{eq:dLSI}
		\ent_\lambda^\delta(\rho)
		\;\leq\;
		\frac{1}{2\alpha^\delta_{\lambda}}
		\delta\sum_k\mob^\delta_k(\rho)\bigl[\dff_k(\ent_\lambda^\delta)'(\rho)\bigr]^2 ,
	\end{equation}
	with $\alpha^\delta_\lambda$ given in~\eqref{eq:alphaLSI}.
\end{proposition}
\begin{proof}
  We establish the modified LSI~\eqref{eq:dLSI} by verifying the curvature condition from~\cite{CaputoPraPosta2009}.
  In this proof, we adopt some notation from~\cite{CaputoPraPosta2009}:
  in particular, we work with the \emph{relative} density $u_\kappa := \rho_\kappa/\bar\rho^\delta_\kappa$
  rather than with $\rho_\kappa$ itself;
  note that
  \[ (\ent_\lambda^\delta)'(\rho)_\kappa = \log u_\kappa
    \quad\text{and}\quad
    \dff_k(\ent_\lambda^\delta)'(\rho) = \dff_k\log u. \]
  In terms of $u$, the relative entropy functional attains the form
  \begin{align*}
    \ent_\lambda^\delta(\rho)
    = \delta\sum_{|\kappa|<K} \bar\rho^\delta_\kappa u_\kappa\log u_\kappa 
    =: \Ent^\delta_{\bar \rho^\delta}(u).
  \end{align*}
  Next, we implicitly introduce a Dirichlet form $\mathcal E$,
  to rewrite the sum on the right-hand side of~\eqref{eq:dLSI}
  as follows:
  \begin{equation}\label{eq:transportId}
    \delta\sum_{|k|<K}\mob^\delta_k(\rho)\bigl[\dff_k(\ent_\lambda^\delta)'(\rho)\bigr]^2
    = -\delta\sum_{|\kappa|<K} \bar\rho^\delta_\kappa \log u_\kappa \mathcal L_\kappa u
    = \mathcal{E}(u,\log u),
  \end{equation}
  where $\mathcal L$ is a (yet to be defined) self-adjoint Markov generator.
  Un-doing the summation by parts that led to the right-hand side of~\eqref{eq:dLSI},
  we obtain
  \begin{align*}
    \delta\mkern-6mu\sum_{|k|<K}\mkern-4mu \mob^\delta_k(\rho)\bigl[\dff_k(\ent_\lambda^\delta)'(\rho)\bigr]^2
    &= -\delta\mkern-6mu\sum_{|k|<K}\mkern-4mu \dff_k\log u\, J_k[\bar\rho^\delta u]
    \\
    &= \delta\mkern-6mu\sum_{|\kappa|<K}\mkern-4mu \bar\rho^\delta_\kappa\,\log u_\kappa\,\frac1{\bar\rho^\delta_\kappa}\dff^*_\kappa J[\rho].
  \end{align*}
  Recalling the definition of $J$ in~\eqref{eq:def:transpFlux}, we see that~\eqref{eq:transportId} is satisfied for all $u:\kgrid^*\to\Rnn$ with the choice
  \begin{align*}
    \mathcal L_\kappa u
    = -\frac1{\bar\rho^\delta_\kappa}\dff^*_\kappa J[\rho]
    = \frac1{\bar\rho^\delta_\kappa}\dff^*_\kappa\bigl(\geom[\bar\rho^\delta]\dff u\bigr),
  \end{align*}
  with the understanding that $\dff_{\pm K}u=0$ because of the boundary condition~\eqref{eq:newbc}.
  It is easily verified that this operator is self-adjoint in the sense that
  \begin{align*}
    \delta\sum_{|\kappa|<K}\bar\rho^\delta_\kappa u_\kappa \mathcal L_\kappa v
    =
    \delta\sum_{|\kappa|<K}\bar\rho^\delta_\kappa v_\kappa \mathcal L_\kappa u
  \end{align*}
  for all functions $u,v:\kgrid^*\to\R$.
  Substituting the explicit expression for~$\bar\rho^\delta$,
  $\mathcal L$ attains the form of a generator for a birth--death chain on $\kgrid^*$,
  \begin{align*}
    \mathcal L_\kappa u = a_\kappa (u_{\kappa+1}-u_\kappa) + b_\kappa (u_{\kappa-1}-u_{\kappa}),
  \end{align*}
  with birth rates $a$ and death rates $b$ given by
  \begin{align*}
    \delta^2a_\kappa
    &= \frac{\geom_{\kappa+1/2}[\bar\rho^\delta]}{\bar\rho^\delta_\kappa}
    = \sqrt{\frac{\bar\rho^\delta_{\kappa+1}}{\bar\rho^\delta_\kappa}}
    = e^{-\frac\lambda 2 \delta x_{\kappa+1/2}}
    \qquad\text{and}\\
    \delta^2b_\kappa
    &= \frac{\geom_{\kappa-1/2}[\bar\rho^\delta]}{\bar\rho^\delta_\kappa}
    = \sqrt{\frac{\bar\rho^\delta_{\kappa-1}}{\bar\rho^\delta_\kappa}}
    = e^{\frac\lambda2 \delta  x_{\kappa-1/2}},
  \end{align*}
  and the understanding that $a_{K-1/2}=0$ and $b_{-K+1/2}=0$;
  the latter is a realization of the boundary conditions $\dff_Ku=0$ and $\dff_{-K}u=0$, respectively.

  We are now in the position to invoke \cite[Theorem 3.1]{CaputoPraPosta2009};
  the result there is stated for one-sided infinite birth--death chains,
  but the result applies to our situation by setting the jump rates to zero outside $\kgrid^*$.
  In terms of the quantities above, the result reads as follows:
  \begin{quotation}
    Assume that $a_\kappa$ and $b_\kappa$ are decreasing and increasing, respectively, with respect to $\kappa\in\kgrid^*$,
    with $a_{K-1/2}=b_{-K+1/2}=0$, and that $c>0$ is such that
    \begin{align}
      \label{eq:caputo}
      a_{k-1/2}-a_{k+1/2} + b_{k+1/2}-b_{k-1/2} \ge c \quad\text{for all $k\in \llbracket -K+1,K-1\rrbracket$}.
    \end{align}
    Then the following modified LSI holds for all functions $u:\kgrid^*\to\Rnn$:
    \begin{align*}
      \Ent^\delta_{\bar \rho^\delta}(u) \le \frac1{2c} \mathcal E(u,\log u).
    \end{align*}
  \end{quotation}
  Monotonicity of $a$ and $b$ is immediately verified.
  The difference in~\eqref{eq:caputo} becomes
  \begin{align*}
    a_{k-1/2}-a_{k+1/2} + b_{k+1/2}-b_{k-1/2}
    = \frac2{\delta^2}\cosh\Bigl(\frac\lambda 2\delta x_k\Bigr)\bigl(1-e^{-\frac\lambda2\delta^2}\bigr).
  \end{align*}
  The expression is minimal for $k=0$, which is an admissible choice and hence we conclude the stated bound~\eqref{eq:alphaLSI}.
\end{proof}
\begin{remark}[Relation to entropic curvature]
	An immediate related question is whether the mixed gradient flow structure based on diffusive and Otto--Wasserstein transport has a suitable notion of $\lambda$-convexity. Note that for the unrescaled DLSS equation on the full space, the results in~\cite{MatthesRottSchlichting2025} suggest that the Otto--Wasserstein metric part is crucial to obtain convexity (see also the numerical experiments in Section~\ref{sec:AnaNumExamples}).

	Let us point out that the discrete relative entropy is indeed geodesically convex with respect to the discrete Otto--Wasserstein metric induced by the transport part of~\eqref{eq:scheme}, see~\cite{ErbarMaas2012,Mielke2013}.
	In particular, an application of~\cite[Theorem 5.1]{Mielke2013} or~\cite[Theorem 4.1]{FathiMaas2016} gives exactly the same lower bound on $\alpha_{\lambda}^{\delta}$ as in~\eqref{eq:alphaLSI} for $\delta>0$. 
  
  Indeed, using the notation from the proof of Proposition~\ref{prop:dLSI}, we can directly apply~\cite[Theorem 5.1]{Mielke2013} or~\cite[Theorem 4.1]{FathiMaas2016} to obtain that the corresponding Markov chain has for any $\lambda>0$ and $\delta>0$ Ricci curvature bounded from below by $\alpha_{\lambda}^\delta$. Indeed, the lower bound is given in terms of 
	\[
		\tfrac{1}{2}\bigl(a_{k-1/2}-a_{k+1/2} + b_{k+1/2}-b_{k-1/2}\bigr) + \tfrac{1}{2} \Xi\bigl(a_{k-1/2}-a_{k+1/2} , b_{k+1/2}-b_{k-1/2}\bigr) \,,
	\]
	with $\Xi(\alpha,\beta)= \inf\set*{ \logmean(r,s)\bigl(\tfrac{\alpha}{r}+\tfrac{\beta}{s}\bigr) : r,s>0}$. By monotonicity of the difference $a_{k-1/2}-a_{k+1/2}$ and $b_{k+1/2}-b_{k-1/2}$ the minimum is attained at $k=0$, in which case $a_{-1/2}-a_{1/2}=b_{1/2}-b_{-1/2} = \delta^{-2} \bigl(1-e^{-\tfrac{\lambda}2 \delta^2}\bigr)$. The conclusion follows by observing that $\Xi(\alpha,\alpha)=2\alpha$ for any $\alpha>0$. 	
\end{remark}

\begin{proof}[Proof of Result C]
  Let $\rho^*:[0,\infty)\to\dprb$ be the solution to~\eqref{eq:dCE2}
  obtained in the proof of Result B, see Section~\ref{sct:proofB}.
  Accordingly, $\rho^*$ is the limit of (a subsequence of) solutions $\rho^K:[0,\infty)\to\dprbK$
  to the problem on the finite grid $\delta\kgrid$,
  with initial data $\rho^{K,0}$ satisfying~\eqref{eq:Hlimit}.
  By Lemma~\ref{lem:EntropyDissipation}, using the dissipations~\eqref{eq:dissipation:diff}--\eqref{eq:dissipation:transp}, we have that
  \begin{align*}
    -\frac{\dd}{\dn t}\ent_\lambda^\delta(\rho^K)
    &= \mathcal{D}^\delta_{\mathrm{diff}}(\rho^K) + \lambda \mathcal{D}^\delta_{\mathrm{transp}}(\rho^K) \\
    &\ge
    0 + \lambda\,\delta\sum_{|k|\le K-1} \mob^\delta_k(\rho^K)\Bigl|\dff_k \bigl(\ent_\lambda^\delta\bigr)'(\rho^K)\Bigr|^2,
  \end{align*}
  so the functional inequality~\eqref{eq:dLSI} above yields that
  \begin{align*}
    -\frac{\dd}{\dn t}\ent_\lambda^\delta(\rho^K)
    \ge 2\lambda\alpha^\delta_\lambda \ent_\lambda^\delta(\rho^K).
  \end{align*}
  So, by Gronwall's inequality,
  \begin{align}
    \label{eq:HKdecay0}
    \ent_\lambda^\delta\bigl(\rho^K(t)\bigr) \le e^{-2\lambda\alpha^\delta_\lambda t}\ent_\lambda^\delta\bigl(\rho^{K,0}\bigr).
  \end{align}
  Now we pass to the limit $K\to\infty$ in~\eqref{eq:HKdecay0}.
  More precisely, we consider a subsequence $K_n\to\infty$
  such that $\rho^{K_n}_\kappa(t)\to\rho^*_\kappa(t)$ locally uniformly in $t\ge0$, for each $\kappa\in\Z^*$.
  Using lower semi-continuity of $\ent_\lambda^\delta$, see Lemma~\ref{lem:Hlsc} on the left-hand side,
  and the initial approximation~\eqref{eq:Hlimit} on the right-hand side,
  one obtains the desired estimate~\eqref{eq:dHdecay}.
\end{proof}

\begin{acknowledgement}
This paper is dedicated to \GS{} on the occasion of his sixtieth birthday.
His work on gradient flows, optimal transport, and dissipative evolution
equations has been a lasting source of inspiration for both authors, and
the results presented here build directly on our joint work with him~\cite{MRSS}.
We are grateful for his generosity as a collaborator and for many stimulating
discussions over the years.

The authors thank José Carrillo and Artur Stephan for raising the question about linearization leading to Remark~\ref{rem:linearization} and Delio Mugnolo for providing the reference~\cite{AGRT}.
\end{acknowledgement}

\bibliographystyle{spmpsci}
\bibliography{GiuseppeBirthday}

\end{document}